\documentclass[leqno,a4paper,titlepage,twoside,11pt]{article}

\usepackage{bbm,pifont,eufrak,eucal,mathrsfs,latexsym,graphicx,epsfig}

\usepackage[all]{xy}

\newcounter{infra}[page]

\newenvironment{dem}[1][]{%
{\bf D\'emonstration #1 : }}{%
\hspace*{\fill}\nolinebreak[1]\hspace*{\fill}\underline{\bf Q.e.d.}\\}

\newenvironment{dem*}[1][]{%
{\bf D\'emonstration #1 : }}{%
}

\newenvironment{eq*}{\begin{eqnarray*}}{\end{eqnarray*}}

\newenvironment{liste}{\begin{itemize}}{\end{itemize}}

\newenvironment{rem}{{\bf{Remarque : }}}{}
\newenvironment{rems}[1][]{{\bf{Remarques #1 :}}}{}

\newtheorem{defi}{D\'efinition}

\newtheorem{thm}{Th\'eor\`eme}[section]
\newtheorem{lem}[thm]{Lemme}
\newtheorem{pro}[thm]{Proposition}

\newcommand{\croi}{\times}

\newcommand{\adh}{\overline}

\newcommand{\boufl}{\xymatrix{
\!\! \ar@(ur,dr)
}}

\newcommand{\C}{\ens C}

\newcommand{\cad}{\mbox{\it c-\`a-d }}
\newcommand{\call}{\mathscr}

\newcommand{\cd}{\rangle}
\newcommand{\cf}{{\it cf. }}
\newcommand{\cg}{\langle}

\newcommand{\ch}{\vee}

\newcommand{\codim}{\mathrm{codim}}

\newcommand{\con}{\supseteq}

\newcommand{\cten}{\Box  \hspace*{-1.77ex}\raisebox{.19ex}{$\croi$}}

\newcommand{\diag}{\mathrm{diag}}

 \newcommand{\Dij}{\bigsqcup}
\newcommand{\donne}{\mapsto}

\newcommand{\End}{\mathrm{End}}
\newcommand{\ens}{\mathbbm}

\newcommand{\equi}{\Leftrightarrow}
\newcommand{\et}{\mbox{ et }}
\newcommand{\etc}{{\it etc}}

\newcommand{\exist}{\exists\:}

\newcommand{\Ext}{\soul{\mathrm{Ext}}}

\newcommand{\fstr}{\mathscr{O}}

\newcommand{\G}{{\mathbf X}}

\newcommand{\GXGT}{\tilda{G} \croi \tilda{G}}

\newcommand{\goth}{\mathfrak}

\newcommand{\hp}{\cdot}

\newcommand{\ie}{\emph{i.e. }}

\newcommand{\impliq}{\Rightarrow}

\newcommand{\inc}{\hookrightarrow}

\newcommand{\infra}[1]{\stepcounter{infra}\setcounter{footnote}{\value{infra}}%
\footnote{#1}}

\renewcommand{\int}{\mathrm{int\:}}

\newcommand{\inter}{\bigcap}
\newcommand{\inv}{^{-1}}

\newcommand{\iso}{\simeq}

\newcommand{\lf}{\longrightarrow}

\newcommand{\limd}[1]{\lim_{\sta{\lf}{#1}}}

\newcommand{\limz}{\ma \lim_{ a \to 0}}

\newcommand{\ma}{\displaystyle}

\newcommand{\moins}{\:\setminus\:}

\newcommand{\ort}{\bot}
\newcommand{\Ou}{\mbox{\Ou}}

\newcommand{\p}{\:\:.}
\newcommand{\PP}{\mathbbm{P}}

\newcommand{\Pic}{\mathrm{Pic}}

\newcommand{\plus}{\oplus}
\newcommand{\Plus}{\bigoplus}

\newcommand{\qed}{\hspace*{\fill}\nolinebreak[1]\hspace*{\fill}\underline{\bf Q.e.d.}\\}
\newcommand{\qq}{\forall\:}

\newcommand{\R}{\ens R}
\newcommand{\res}[1]{{\left | {}_{#1} \right.}}

\newcommand{\re}[1]{\textbf{\textsc{#1}}}

\newcommand{\si}{\mathrm{\; si \;}}
\newcommand{\sinon}{\mbox{ sinon }}
\newcommand{\SL}{\mathrm{SL}}

\newcommand{\soul}{\underline}

\newcommand{\sta}{\stackrel}
\newcommand{\sub}{\subseteq}

\newcommand{\surj}{\:\rule[.5ex]{1.5em}{0.1ex}\!\!\!\!\!\succ\!\!\!\succ }
 
\newcommand{\sus}{\ldots}

\newcommand{\tens}{\otimes}
\newcommand{\tenso}[1]{\raisebox{-1.5ex}{$\ma \stackrel{\displaystyle \tens}{\scriptstyle #1}$}}
\newcommand{\ti}{\tilda}
\newcommand{\tilda}{\widetilde}

\newcommand{\tq}{\: : \:}

\newcommand{\tri}{$$
\ast\;\ast\;\ast
$$}

\newcommand{\Uni}{\bigcup}
\newcommand{\uni}{\bigcup}
\newcommand{\unio}{\cup}

\newcommand{\Vect}{\bigwedge\nolimits}

\newcommand{\vide}{\emptyset}

\newcommand{\X}{{\cal X}}

\newcommand{\Y}{{\cal Y}_\R}

\newcommand{\Z}{\ens Z}
\newcommand{\z}{\mbox{\bf z}}

\newcommand{\appli}[5][]{\begin{array}{ccc}
#2 & \stackrel{#1}{\longrightarrow} & #3 \\
#4 & \longmapsto & #5 
\end{array}}

\newcommand{\applinom}[6][]{#2:\left\{
\begin{array}{ccc}
#3 & \stackrel{#1}{\longrightarrow} & #4 \\
#5 & \longmapsto & #6
\end{array}
\right.
}

\begin{document}

\title{Cohomologie des fibr\'es en droites sur les compactifications des groupes r\'eductifs}

\author{Alexis TCHOUDJEM \\
Institut Fourier\\
Universit\'e de Grenoble\\
100, rue des maths\\
B.P. 74\\
38402 Saint-Martin d'H\`eres cedex\\
atchoudj@ujf-grenoble.fr}

\maketitle

\tableofcontents

\section*{Introduction}

Pour un groupe $G$ r\'eductif et connexe sur $\C$, le c\'el\`ebre th\'eor\`eme de Borel-Weil-Bott d\'ecrit tr\`es simplement les groupes de cohomologie des fibr\'es en droites sur les $G-$vari\'et\'es de drapeaux, au moyen des repr\'esentations irr\'eductibles de $G$. Notre but est d'\'etablir un r\'esultat analogue pour une autre famille de vari\'et\'es : les compactifications \'equivariantes du groupe $G$. Comme celui-ci est un espace homog\`ene pour l'action, par multiplication \`a gauche et \`a droite, de $G \croi G$, ces compactifications admettent une description combinatoire qui r\'esulte de la th\'eorie des plongements des espaces homog\`enes sph\'eriques (\cf \cite{LV:peh}). On peut m\^eme obtenir une description combinatoire des fibr\'es en droites $L$ sur une compactification \'equivariante $X$.

Nous exprimons ici les groupes de cohomologie de $L$ sur $X$, en fonction de ces donn\'ees combinatoires. Pour un rev\^etement fini $\ti{G}$ de $G$, tel que $L$ est $\GXGT-$lin\'earisable\infra{
Un tel rev\^etement existe toujours \cf \cite[remarque p.67]{KKLV}.}, ces groupes $H^i(X,L)$ sont en fait des $\GXGT-$modules (pour tout $i \ge 0$). Notre th\'eor\`eme principal (\ref{thm:princi}) donne leurs multiplicit\'es selon chaque $\GXGT-$module simple. Ce r\'esultat fait intervenir des $\ll$ objets toriques $\gg$ et, pour y parvenir, j'utilise le complexe de Grothendiek-Cousin, tel qu'il a \'et\'e introduit notamment dans \cite{Ke:GCc}. C'est un complexe fait de groupes de cohomologie \`a support dans les $B \croi B^-$ orbites de $X$ (on note $B$ et $B^-$ deux sous-groupes de Borel oppos\'es de $G$) et dont l'homologie est exactement la cohomologie $H^*(X,L)$.

Dans le contexte des vari\'et\'es de drapeaux, le complexe correspondant permet de retrouver directement le th\'eor\`eme de Borel-Weil-Bott. Pour ce qui concerne les compactifications, on a besoin d'analyser un peu plus les groupes de cohomologie \`a support qui apparaissent. Si $\goth g$ est l'alg\`ebre de Lie de $G$, ce sont des repr\'esentations de $\goth g \croi \goth g$. On en donnera des filtrations avec des gradu\'es associ\'es plus familiers, en se servant de la structure de la compactification. Cela suffira pour calculer les multiplicit\'es cherch\'ees.

Auparavant (\cf la proposition \ref{pro:ep}), ces filtrations auront mis en lumi\`ere l'utilit\'e de consid\'erer $G$ dans une compactification, pour pr\'eciser la structure de certains $\goth g \croi \goth g-$modules int\'eressants du point de vue de la th\'eorie des repr\'esentations, comme, par exemple, les groupes de cohomologie de $G$ \`a support dans les doubles classes ${B w B^-}$.

Lorsque $G$ est un tore, les compactifications sont les vari\'et\'es toriques compl\`etes ; dans ce cas, les groupes de cohomologie des fibr\'es en droites ont \'et\'e d\'etermin\'es par Demazure (\cf \cite{De:sg}). Notre r\'esultat est aussi une g\'en\'eralisation de ce th\'eor\`eme de Demazure.

D'autre part, lorsque $G$ est semi-simple adjoint, une compactification particuli\`ere de $G$ est sa compactification magnifique construite par De Concini et Procesi. Pour cette compactification, notre r\'esultat a \'et\'e d\'emontr\'e par la m\^eme m\'ethode, de mani\`ere ind\'ependante par S. Kato (\cf \cite{K:BWB}), et aussi annonc\'e dans \cite{moi}.

\section*{Notations}

{\it  Ici \label{par:ici} $G$ sera un groupe r\'eductif et connexe sur $\C$, d'alg\`ebre de Lie $\goth g$. On choisira $B$ un de ses sous-groupes de Borel et $T$ un tore maximal de $B$ ; on appellera $B^-$ le sous-groupe de Borel oppos\'e \`a $B$, relativement \`a $T$ (\ie tel que $B^- \cap B = T$). Soient $\Phi$ et $W$ le syst\`eme de racines et le groupe de Weyl de $(G,T)$. On notera $\Phi^+$  l'ensemble des racines positives, $\Phi^-$ celui des racines n\'egatives, $\rho$ la demi-somme des racines positives, et, si $w \in W$, $w * \lambda := w(\lambda + \rho) - \rho$ pour tout caract\`ere $\lambda$ de $T$. Soient $\Delta$ la base de $\Phi$ d\'efinie par $B$, $l$ la fonction longueur correspondante sur $W$ et $w_0$ l'\'el\'ement le plus long de $W$. Le rang de $G$, \cad le cardinal de $\Delta$, sera not\'e $r$.}

\begin{center}
\it
Les vari\'et\'es consid\'er\'ees dans ce texte seront des vari\'et\'es alg\'ebriques sur $\C$.
 \end{center}

\section{Rappels sur les compactifications de groupes}

Apr\`es la d\'efinition g\'en\'erale des compatifications de groupes, on s'int\'eressera surtout aux compactifications r\'eguli\`eres. On rappelle cette notion en \ref{sec:cr} et aussi les donn\'ees combinatoires qui caract\'erisent les faisceaux inversibles sur ces vari\'et\'es.

\subsection{D\'efinition g\'en\'erale}\label{sec:dg}

{\it
On appelle compactification (\'equivariante) de $G$ toute vari\'et\'e compl\`ete, normale qui contient $G$ comme ouvert et o\`u l'action par multiplication \`a gauche et \`a droite de $G \croi G$ sur $G$ se prolonge.}

\tri

Par exemple, si   $G = \SL(2,\C)$ et si $M(2,\C)$ est l'espace des matrices d'ordre $2$, la quadrique de $\PP^4$ :
$${\cal Q} := \left\{ [M:t] \in \PP(M(2,\C)) \croi \C \tq \det M = t^2 \right\}$$
est munie d'une action de $G \croi G$ d\'efinie par :
$$
\appli{ G\croi G \croi {\cal Q} }{ \cal Q }{ (g_0,g_1,[M:t]) }{ [g_0 M g_1\inv : t ] } \p$$
Pour cette action et pour l'inclusion :
$$\appli{ G } {\cal Q}{ M } { [M : 1] } \;\; , $$
${\cal Q}$ est une compactification lisse de $G$.

\tri

Voici un autre exemple : si $G = PGL(2, \C) = SL(2,\C) / \{\pm Id\}$, alors l'espace projectif $\PP (M(2, \C))$ est une compactification lisse de $G$.

\tri

D'apr\`es \cite{B:Inf}, dans une compactification du groupe $G$, tout point fix\'e par le tore maximal $T \croi T$ appartient \`a une orbite ferm\'ee de $G \croi G$. Nous verrons plus loin (\cf la partie \ref{sec:sdc}) que c'est un avantage, par rapport aux compactifications d'autres espaces homog\`enes sph\'eriques, pour d\'eterminer la cohomologie des fibr\'es en droites.

\subsection{Compactifications r\'eguli\`eres }\label{sec:cr}

La notion de vari\'et\'e r\'eguli\`ere a \'et\'e introduite par Bifet, De Concini et Procesi dans \cite{BDP:cre}, \cf aussi \cite[\S 1.4]{B:Inf} et \cite{BB:alrrv} :

\begin{defi}\label{def:vreg}
On dit qu'une vari\'et\'e $X$, munie d'une action de $G$ est \re{r\'eguli\`ere} si sont v\'erifi\'ees :
\begin{enumerate} 
\item $X$ est lisse et a une $G-$orbite ouverte et dense $X^0_G$ dont le compl\'ementaire est un {\it diviseur \`a croisements normaux}. On appellera  \re{diviseurs limitrophes} les composantes irr\'eductibles de  $X - X^0_G$.

\item Chaque adh\'erence de $G-$orbite est l'{\it intersection transverse} des diviseurs limitrophes qui la contiennent.

\item Si $x \in X$ alors sur l'espace normal \`a $G \cdot x$ dans $X$, le groupe d'isotropie de $x$, agit avec une orbite dense.
\end{enumerate}
\end{defi}

\begin{defi}
Les compactifications r\'eguli\`eres de $G$ sont les compactifications (\'equivariantes) $X$ de $G$ qui sont r\'eguli\`eres comme $G \croi G-$vari\'et\'es.
 \end{defi}

{\bf Exemples : }

--- Les vari\'et\'es toriques compl\`etes lisses sont les compactifications r\'eguli\`eres des tores.

--- La compactification ${\cal Q}$ de \ref{sec:dg} est une compactification r\'eguli\`ere.

--- Pour un groupe adjoint $G$ (\ie de centre trivial), De Concini et Procesi ont construit (Cf. \cite{DeCP:CSV}) une compactification lisse, $\ll$ canonique $\gg$, ou $\ll$ magnifique $\gg$, $\adh{G}$, de $G$ ; on peut la d\'efinir comme l'unique compactification r\'eguli\`ere de $G$ avec une seule $G \croi G-$orbite ferm\'ee. L'espace projectif $\PP(M(2,\C))$ est, par exemple, la compactification magnifique de $PGL(2,\C)$.

\tri

\begin{center}
\sc Jusqu\`a la fin $\G$ sera une compactification r\'eguli\`ere de $G$.
 \end{center}

Avant d'\'enoncer le th\'eor\`eme principal sur les groupes de cohomologie, fixons les notations qui concernent $\G$ :

\subsubsection{Caract\'erisation combinatoire de la compactification}\label{sec:dc}

On fait d'abord appel \`a deux r\'eseaux en dualit\'e : le r\'eseau ${\cal X}$ des caract\`eres de $T$ et  celui, ${\cal Y}$, des sous-groupes \`a un param\`etre de $T$. On posera $\ma {\cal X}_\R := {\cal X} \tenso{\Z} \R$ et $\ma {\cal Y}_\R := {\cal Y} \tenso{\Z} \R$. On notera :
$$\cg \: , \: \cd : {\cal X}_\R \croi {\cal Y}_\R \to \R$$
le crochet de dualit\'e usuel.

 Soit  :
$${\cal C}^+ := \left\{ \nu \in \Y \tq \qq \alpha \in \Delta , \, \cg \alpha , \nu \cd \ge 0 \right\}$$
l'adh\'erence de la chambre de Weyl positive.

\tri

On va d\'efinir une subdivision de ${{\cal C}^+}$ associ\'ee \`a  $\G$.

Pour cela, soit $\adh{T}$ l'adh\'erence de $T$ dans $\G$. Sur $G$, la restriction de l'action de $G \croi G$ \`a la diagonale de $T$, $\diag (T)$,  est donn\'ee par :
$$\qq g \in G, \: \qq t \in T, \; (t,t). g = tgt\inv \p$$
Elle se prolonge \`a  $\G$. La vari\'et\'e des points fix\'es par le tore diagonal $\diag(T)$ est lisse (\cf \cite[prop. 1.3]{I:fpf}) et $\adh{T}$ en est une composante irr\'eductible.  Pour l'action de $T$ \`a gauche (\ie pour l'action de $T \croi \{1\}$),  $\adh{T}$ est donc une vari\'et\'e torique compl\`ete lisse  ; on note ${\cal E}$ {\it l'\'eventail associ\'e \`a $\adh{T}$}. Comme $\adh{T}$  est compl\`ete, ${\cal E}$ est une subdivision de ${\cal Y}_\R$ (\ie : $\ma \uni_{ \sigma \in {\cal E} } \sigma = {\cal Y}_\R$) et comme $\adh{T}$ est invariant par l'action de la diagonale du  groupe de Weyl :
 $$\diag (W) := \{(w,w) \tq w \in W\} \;\;,$$
 ${\cal E}$ est aussi $W-$invariant. Il r\'esulte de \cite[prop. A2]{B:Inf} que ${\cal E} = W {\cal E}^+$ o\`u ${\cal E}^+$ est la subdivision de ${\cal C}^+$ form\'ee des c\^ones de ${\cal E}$ contenus dans ${{\cal C}^+}$. 
De mani\`ere analogue au cas torique, les c\^ones de l'\'eventail ${\cal E}^+$ param\`etrent des orbites : les $G \croi G-$orbites de $\G$ (\cf encore \cite[prop. A2]{B:Inf}). Pour tout $\sigma \in {\cal E}^+$, on d\'esigne par $z_\sigma$ le point-base correspondant dans $\adh{T}$  et par ${\cal O}_\sigma$ sa $G \croi G-$orbite. On dira que $z_\sigma$ est le {\it \label{pb} point-base de l'orbite ${\cal O}_\sigma$}. Parfois on notera $z_{{\cal O}_\sigma}$ le point $z_\sigma$. 

Les $G \croi G-$orbites ferm\'ees de $\G$ correspondent aux c\^ones de ${\cal E}^+$ de dimension maximale $r$. Lorsque $\sigma$ est un tel c\^one, $z_\sigma$ a pour groupe d'isotropie $B^- \croi B$ ; d'o\`u :
$$(G \croi G ) \hp z_\sigma \iso G / B^- \croi G/B $$
(\cf \cite[prop. A1 ii)]{B:Inf}).

\begin{rem}
R\'eciproquement, toute subdivision de ${\cal C}^+$ dont tous les c\^ones sont lisses (\cad engendr\'es par un d\'ebut de base du r\'eseau ${\cal Y}$) correspond \`a une compactification r\'eguli\`ere de $G$. Par exemple, la subdivision triviale form\'ee de ${\cal C}^+$ et de ses faces d\'efinit la compactification magnifique de $G$ (si $G$ est adjoint).
 \end{rem}

\tri

On va maintenant pr\'eciser les $\ll$ objets toriques $\gg$ correspondant aux faisceaux inversibles sur $\G$.

\subsubsection{Caract\'erisation combinatoire des faisceaux inversibles}\label{sec:fi}

On choisit au pr\'ealable une isog\'enie $r : \ti{G} \to G$ telle que $\Pic(\ti{G}) = (0)$ (\cf \cite{I:GAG} ou \cite[prop. 4.6]{KKLV}). Tous les faisceaux inversibles sur $\G$ sont alors $\GXGT-$lin\'earisables.

On notera $\ti{B}$, $\ti{B^-}$, $\ti{T}$, $\ti{\cal X}$, $\ti{\cal Y}$ les ensembles correspondants \`a ${B}$, ${B^-}$, ${T}$, ${\cal X}$, ${\cal Y}$ pour $\ti{G}$. 

Gr\^ace \`a la param\'etrisation des fibr\'es en droites sur les vari\'et\'es sph\'eriques de \cite[\S 2.2]{B:nc}, on obtient\infra{
On utilise en particulier que $\G$ admet un recouvrement par des ouverts isomorphes \`a des espaces affines et les trivialisations des fibr\'es en droites restreints \`a ces ouverts. 
} celle des fibr\'es en droites sur $\G$ par les {\it fonctions lin\'eaires par morceaux sur ${{\cal C}^+}$, adapt\'ees \`a la subdivision ${\cal E}^+$}. Il s'agit des applications $h : {\cal C}^+ \to \R$ telles que pour tout c\^one $\sigma$ de ${\cal E}^+$ et tout $n \in \sigma$, $h(n) = \cg h_\sigma , n \cd$, pour un certain $\ma h_\sigma \in {\cal X}_\R$.

\begin{thm}[{ \cite[\S 2.2]{B:nc} }]\label{thm:parfdc}

\'Etant donn\'ee une fonction lin\'eaire par morceaux  $ h : {{\cal C}^+} \to \R$ \`a valeurs enti\`eres sur $\ti{\cal Y} \cap {\cal C}^+$ et adapt\'ee \`a l'\'eventail ${\cal E}^+$, il existe un  faisceau $\call{L}$, inversible et $\ti{G} \croi \ti{G}-$lin\'earis\'e sur $\G$, tel que, pour toute orbite ferm\'ee ${\cal O}_\sigma$ ($\sigma$ c\^one maximal de ${\cal E}^+$), le groupe $\ti{B^-} \croi \ti{B}$ op\`ere avec le caract\`ere $(h_\sigma, -h_\sigma)$ dans la fibre $\call{L}\res{z_\sigma}$.

A isomorphisme de faisceaux $\ti{G} \croi \ti{G}-$lin\'earis\'es sur $\G$ pr\`es, $\call{L}$ est unique et on le note : $\call{L}_h$.

On obtient ainsi, \`a isomorphisme de $\fstr_{\G}-$modules pr\`es, tous les faisceaux inversibles sur $\G$ $\Box$

 \end{thm}

{\sc On conservera cette notation $\call{L}_h$ pour le faisceau inversible et $\GXGT-$lin\'earis\'e sur $\G$ correspondant \`a la fonction lin\'eaire par morceaux $h$.}

\begin{rem}
D'apr\`es \cite[lemme 2.2]{KKV}, deux lin\'earisations (pour $\GXGT$) d'un m\^eme faisceau inversible $\call{L}$ sur $\G$ diff\`erent d'un caract\`ere de $\GXGT$. En cons\'equence, c'est \`a translation pr\`es par un caract\`ere de $\ti{G}$ que le faisceau $\call{L}$ d\'efinit la fonction $h$.
 \end{rem}

\section{Description de la cohomologie des faisceaux inversibles sur une compactification de groupe}

En premier lieu, on \'enonce le th\'eor\`eme principal pour une compactification r\'eguli\`ere. Ensuite, on en donnera une version dans le cas particulier des compactifications magnifiques. 

\begin{rem}
Si $X$ est une compactification \'equivariante quelconque de $G$ (\ie seulement normale), alors, d'apr\`es \cite{B:ext}, il existe une compactification r\'eguli\`ere $\ti{\G}$ de $G$ et un morphisme propre et birationnel $\pi : \ti{\G} \to X$ tel que :
$$1) \; \pi_*\fstr_{\ti{\G}} = \fstr_X \;\;\;\; 2) \;\qq i > 0 , \; R^i\pi_* \fstr_{\ti{\G}} = (0) \p$$ En particulier, pour tout faisceau inversible $\call{L}$ sur $X$, on a : $H^i(X, \call{L}) = H^i(\ti{\G} , \pi^*\call{L} )$ pour tout $i$. Notre th\'eor\`eme \ref{thm:princi} est donc valable, en toute g\'en\'eralit\'e, pour les groupes de cohomologie des fibr\'es en droites sur les compactifications \'equivariantes de $G$.
 \end{rem}

\subsection{Cas r\'egulier}

Soit $h$ une fonction lin\'eaire par morceaux sur ${{\cal C}^+}$, adapt\'ee \`a ${\cal E}^+$ et \`a valeurs enti\`eres sur $\ti{\cal Y} \cap {\cal C}^+$.

Les groupes de cohomologie $H^i(\G, \call{L}_h)$ sont des modules de dimension finie sous le  groupe $\GXGT$. Ils se d\'ecomposent donc en somme directe de repr\'esentations irr\'eductibles. Notons ${\cal X}^+$ (respectivement $\ti{\cal X}^+$) l'ensemble des caract\`eres dominants de ${\cal X}$ (respectivement de $\ti{\cal X}$). Pour chaque $\mu \in \ti{\cal X}^+$, soit $L(\mu)$ le $\ti{G}-$module simple de plus haut poids $\mu$. Les seuls modules simples apparaissant dans la d\'ecomposition des groupes de cohomologie $H^i(\G, \call{L}_h)$ sont les 
$$\End(L(\mu)) = L(\mu) \tens L(-w_0 \mu) ,\;\; \mu \in \ti{\cal X}^+\p$$

Nous allons exprimer leurs multiplicit\'es au moyen des ensembles :
$$ V(h,\lambda) := \{ n \in {{\cal C}^+} \tq \cg \lambda , n \cd - h(n) > 0 \} $$
o\`u $\lambda \in {\cal X}_\R$. Pour abr\'eger, on note $h + {\cal X}$ l'ensemble des caract\`eres $\mu \in {\cal X}_\R$ tels que $\mu - h_\sigma \in {\cal X}$ pour tout c\^one maximal $\sigma$ de ${\cal E}^+$. On pose aussi 
$$J_t : = \{\alpha \in \Delta \tq t\inv(\alpha) < 0 \}$$ et, pour chaque $\alpha \in \Delta$,  $$\alpha^\ort : = \{n \in {{\cal C}^+} \tq \cg \alpha , n \cd = 0 \} \p$$ 

On utilisera  la notation ${H}^*( \; , \; )$ pour les groupes de cohomologie relative de parties de ${\cal Y_\R}$ \`a valeurs dans $\R$ (\cf par exemple \cite[cha. 5, sec. 4, ex.5]{Spa:at}).  

Cela \'etant pos\'e :

\begin{thm}\label{thm:princi}
On a un isomorphisme de $\GXGT-$modules :
$$H^i(\G , \call{L}_h) \iso \Plus_{\mu \in \ti{\X}^+} m^i_h(\mu) \End(L(\mu))$$
o\`u pour tout poids entier  dominant $\mu$ et tout entier $i$, la multiplicit\'e $m^i_h(\mu)$ est \'egale \`a :

$$\sum_{t \in W \moins \{1\}} \dim {H}^{i - 2 l(t)-1} \left( V(h, t*\mu   ) , V( h,t*\mu  ) \cap \uni_{\alpha \in J_t} \alpha^\ort \right)$$
$$  + \dim H^i({{\cal C}^+} , V(h,\mu)) $$
si $\mu \in h + {\cal X}$ et $m^i_h(\mu) = 0$ sinon.
 \end{thm}

\begin{center}
*
 \end{center}

Afin d'illustrer cette formule compliqu\'ee, voici quelques remarques et exemples :

\begin{liste}

\item[1.---] Avec un $\ti{G}$ choisi de sorte que pour tout $\mu \in \ti{X}$ et pour toute coracine $\alpha^\ch$, on ait : $\cg \mu , \alpha^{\ch} \cd \in \Z$, la condition $\mu \in h + {\cal X}$ s'\'ecrit :
$$\qq n \in {\cal Y} \cap {\cal C}^+, \; \cg \mu , n \cd - h(n) \in \Z \p$$ 

 \item[2.---] Lorsque $G = T$, on a dans ce cas :
$$W = \{1\} , \; \Delta = \vide \et {\cal C}^+ = {\cal Y}_\R \p$$
On retrouve le r\'esultat de Demazure (\cf par exemple \cite[th. 2.6]{O:cbag}) :
$$m^i_h(\mu) = \dim H^i({\cal Y}_\R , V(h,\mu) ) \p$$
Toutefois, ce r\'esultat est en fait utilis\'e dans la d\'emonstration (\cf la section \ref{sec:ract}).

\item[3.---] Lorsque $i = 0$, on obtient tout de suite :

$$m^0_h(\mu) = \left\{ \begin{array}{ll}
1 & \mbox{ si pour tout $n \in {\cal C}^+$, $\cg \mu , n \cd \le h(n)$ et si $\mu \in h + {\cal X}$} \\
0 & \sinon.
 \end{array}\right.$$

C'est un r\'esultat d\'ej\`a connu (\cf par exemple \cite[th. du \S 3.4]{Bi:ocsv}).

\item[5.---] Soit $X({\cal E}^+)$ la vari\'et\'e torique associ\'ee \`a l'\'eventail ${\cal E}^+$. C'est un ouvert de $\adh{T}$ et, lorsque $i= 0 , 1,$ ou $2$, on obtient :
$$m^i_h(\mu) = \dim \left(H^i(X({\cal E}^+), \call{L}_h)\right)_\mu$$
la dimension de l'espace propre associ\'e au caract\`ere $\mu$.

\tri

\item[5.---]
On va maintenant voir un exemple de calcul o\`u l'on trouve une multiplicit\'e $>1$. Lorsque $ G = PGL(3 , \C)$, et $\ti{G} = SL(3,\C)$, les groupes de cohomologie des fibr\'es en droites sur la compactification magnifique $\adh{G}$ sont des repr\'esentations de $SL(3,\C) \croi SL(3,\C)$ sans multiplicit\'e (\cad que leurs multiplicit\'es sont $0$ ou $1$), \cf la remarque 4 p. \pageref{rem:magni4}. Ce n'est plus le cas pour les compactifications g\'en\'erales du m\^eme groupe.

  Soient $\alpha,\beta$ les \'el\'ements de la base $\Delta$. On note $\{\omega_\alpha^\ch , \omega_\beta^\ch \}$ la base duale de $\{\alpha,\beta\}$ dans ${\cal Y}$.

Soit $\G$ l'\'eclat\'e de la compactification magnifique de $G$ le long de son unique $G \croi G-$orbite ferm\'ee. C'est aussi la compactification r\'eguli\`ere de $G$ associ\'ee \`a la subdivision suivante de la chambre de Weyl positive :
$$
{{\cal C}^+} = \sigma^\alpha \uni \sigma^\beta\;,$$ 
o\`u :
$$\sigma^\alpha : = \R_+ \omega_\alpha^\ch  + \R_+ (\omega_\alpha^\ch + \omega_\beta^\ch )$$
et :
$$ \sigma^\beta :=  \R_+ (\omega_\alpha^\ch + \omega_\beta^\ch) + \R_+ \omega_\beta^\ch \p$$


\begin{figure}[hbtp]
  \begin{center}
    \leavevmode
    \input{ecla.pstex_t}
    \caption{L'\'eventail ${\cal E^+}$}
  \end{center}
\end{figure}


Soit $h$ la fonction d\'efinie par :

$$h(n) = \left\{ \begin{array}{ll}
\cg \beta - 2\alpha , n \cd & \si n \in \sigma^\alpha \\
\cg \alpha - 2\beta , n \cd & \si n \in \sigma^\beta
\end{array}
\right.
$$(il n'y a pas d'ambigu\"{i}t\'e).
Notons $\call{L}_h$ le faisceau inversible correspondant sur $\G$.

On va montrer que la multiplicit\'e de la repr\'esentation triviale $\End(L(0))$ dans $H^3( \G , \call{L}_h )$ vaut $2$.

On remarque d'abord que le caract\`ere $0 \in h + {\cal X}$. Ensuite, puisque si $t \in W$ :
$$ 3-2l(t) - 1 \ge 0 \equi l(t) \le 1 \;,
$$
on a d'apr\`es le m\^eme th\'eor\`eme :
$$m_h^3(0) = \sum_{t \in \{s_\alpha ,s_\beta\}} \dim H^{0} \left( V(h,t* 0) , V( h , t *0) \cap \uni_{\delta \in J_t} \delta^\ort  \right)$$
$$ + \dim H^3({{\cal C}^+}, V(h,0) ) \p$$
Or, d'une part : 
$$V(h,0) = \{ n \in {{\cal C}^+} \tq h(n) <0 \}$$
$$ = \left\{  n_\alpha \omega_\alpha^\ch + n_\beta \omega_\beta^\ch \tq (n_\alpha,n_\beta) \in \R_+^2 \moins \{(0,0)\} \right\}$$

donc : $H^3({{\cal C}^+} , V(h,0) ) = (0)$.

D'autre part, on montre que :
$$V(h,s_\alpha * 0) = \{ n \in {{\cal C}^+} \tq - \cg \alpha , n \cd - h(n) > 0 \}$$
$$ = \left\{  n_\alpha \omega_\alpha^\ch + n_\beta \omega_\beta^\ch \tq (n_\alpha,n_\beta) \in \R_+^2 \moins \{(x,x) \tq x \ge 0\} \right\}$$ et :
$$V(h,s_\alpha * 0) \cap \alpha^\ort = \R_+^* \omega_\beta\p$$

D'o\`u :
$$H^0(V(h,s_\alpha * 0) , V(h,s_\alpha * 0) \cap \alpha^\ort  ) $$
$$= H^0(\R_+ \croi \R_+^* , \{0\} \croi \R_+^* ) \plus H^0(\R_+^* \croi \R_+ , \vide )  =\R \p$$

On montre de m\^eme, que :
$$V(h,s_\beta * 0) = \left\{  n_\alpha \omega_\alpha^\ch + n_\beta \omega_\beta^\ch \tq  (n_\alpha,n_\beta) \in \R_+^2 \moins \{(x,x) \tq x \ge 0\} \right\}$$
et que :
$$V(h,s_\beta * 0) \cap \beta^\ort = \R_+^* \omega_\alpha \p$$

Donc : $\ma H^0(V(h,s_\beta * 0) , V(h,s_\beta * 0) \cap \beta^\ort ) = \R \p$

On a ainsi : $ \ma m^3_h(0) = 2$.
\qed

\end{liste}

\subsection{Cas de la compactification magnifique}
 Soit $G_{ad}$ le groupe adjoint de $G$ (le quotient de $G$ par son centre) et $\adh{G_{ad}}$ sa compactification magnifique (\cf la section \ref{sec:cr}). On peut prendre pour $\ti{G_{ad}}$ le rev\^etement universel de $G_{ad}$. Les faisceaux inversibles sur $\adh{G_{ad}}$ sont alors simplement donn\'es par un poids entier $\lambda \in \ti{\cal X}$. Soit $z$ l'unique point fixe de $\adh{G_{ad}}$ pour $B^- \croi B$ et pour chaque $\lambda \in \ti{\cal X}$, soit $\call{L}_\lambda$ le faisceau inversible et $\ti{G_{ad}} \croi \ti{G_{ad}}-$lin\'earis\'e sur $\adh{G_{ad}}$ tel que le groupe $\ti{B^-} \croi \ti{B}$ op\`ere dans la fibre $\call{L}\res{z}$ avec le caract\`ere $(\lambda, -\lambda)$.

On note aussi $Q : = \sum_{\alpha \in \Delta} \Z \alpha$ le r\'eseau radiciel et pour chaque $t \in W$, $Q_t := \left\{  \sum_{\alpha \in \Delta} \gamma_\alpha \alpha \in Q \tq \gamma_\alpha > 0 \equi t\inv(\alpha) <0 \right\} \p$ 

On a :
\begin{thm}\label{thm:magni}
Comme $\ti{G_{ad}} \croi \ti{G_{ad}}-$modules, pour tout $i \ge 0$ :
$$H^i(\adh{G_{ad}} , \call{L}_\lambda) = \Plus_{\mu \in \ti{\cal X}^+} m_\lambda^i(\mu) \End(L(\mu))$$
avec chaque multiplicit\'e $m^i_\lambda(\mu)$ valant le nombre de $t\in W$ tels qu'\`a la fois  $2l(t) + |J_t|=i$ et $t * \mu \in \lambda + Q_t$.
 \end{thm}

\tri

Ce th\'eor\`eme \ref{thm:magni} est un cas particulier du th\'eor\`eme \ref{thm:princi} comme nous allons le voir tout de suite. On remarque d'abord que, dans ce cas, $h$ est simplement un caract\`ere $\lambda$ de $\ti{\cal X}$. Soit $\mu \in \ti{\cal X}^+$. On a alors : $$V(\lambda,\mu) =  \{ n \in {{\cal C}^+} \tq \cg \mu , n \cd - \cg \lambda , n \cd > 0 \}\p$$
 La condition $\mu \in h + {\cal X}$ du th\'eor\`eme signifie que $\mu - \lambda \in {\cal X}$.

On fixe ensuite $ \nu \in {\cal X}$. Soit $V := \{ n \in {{\cal C}^+} \tq \cg \nu , n \cd > 0 \}$. On va montrer (et cela suffira) que, pour tout $j \ge 0$ et tout $ t \in W \moins \{1\}$ :
$$ H^j(V , V \cap \uni_{\alpha \in J_t} \alpha^\ort) = \R \si j = |J_t| - 1 \et \si \nu \in Q_t$$
et que $\ma H^j (V , V \cap \uni_{\alpha \in J_t} \alpha^\ort) = (0)$ dans tous les autres cas.

Comme cela est vrai lorsque $ V = \vide$, on va maintenant supposer que $V$ n'est pas vide.

Puisque $V$ est convexe, on a, pour tout $j$ :
$$H^j(V , V \cap \uni_{\alpha \in J_t} \alpha^\ort) = \ti{H}^{j-1} ( V \cap \uni_{\alpha \in J_t} \alpha^\ort)$$
o\`u $\ti{H}^*(\;)$ d\'esigne la cohomologie r\'eduite (\cf par exemple  \cite[Chap. 5, sec 4, \S 2]{Spa:at}).

Avec la base duale $\{ \omega_\alpha^\ch \tq \alpha \in \Delta\}$ de $\Delta$ dans ${\cal Y}$, on d\'efinit une fonction :
$$\ma \applinom{p}{{{\cal C}^+}}{\R_+^{J_t}}{n= \sum_{\delta \in \Delta} n_\delta \omega_\delta^\ch}{\sum_{\delta \in J_t} n_\delta \omega_\delta^\ch} \p$$

Cette application $p$ est continue et on s'aper\c{c}oit que, restreinte \`a $V \cap \uni_{\alpha \in J_t} \alpha^\ort$, ses fibres sont convexes :
$$ \qq n^0 \in \R_+^{ J_t}, \; \{ n \in V \cap \uni_{\alpha \in J_t} \alpha^\ort \tq p(n) = n^0 \}$$
$$ = \vide \mbox{ ou } \{ n = n^1 + n^0 \tq n^1 \in \R_+^{\Delta \moins J_t} \et \cg \nu , n^1 \cd > - \cg \nu , n^0 \cd \} \p$$

On en d\'eduit gr\^ace \`a une suite spectrale de Leray (\cf \cite[chap. II, th\'eor\`eme 4.17.1]{Go}) que :
$$\ti{H}^{j-1} ( V \cap \uni_{\alpha \in J_t} \alpha^\ort ) = \ti{H}^{j-1}(p(V \cap \uni_{\alpha \in J_t} \alpha^\ort) )$$
pour tout $j$.

Mais, si on note $\partial \R_+^{J_t}$ le bord de $\R_+^{J_t}$, sont v\'erifi\'ees :
$$p(V \cap \uni_{\alpha \in J_t} \alpha^\ort) = \{ n \in \R_+^{J_t} \cap \uni_{\alpha \in J_t} \alpha^\ort \tq \exist n^1 \in \R_+^{\Delta \moins J_t} , \; \cg \nu , n \cd + \cg \nu , n^1 \cd > 0 \}$$
$$ = \R_+{J_t} \cap \left( \uni_{\alpha \in J_t} \alpha^\ort \moins \{ n \in \R_+^{J_t} \cap \uni_{\alpha \in J_t} \alpha^\ort \tq \cg \nu ,n \cd \le - \sup_{n^1 \in \R_+^{\Delta \moins J_t} } \cg \nu , n^1 \cd \} \right)$$
$$ = \partial \R_+^{J_t} \mbox{ ou } \partial \R_+^{J_t} \moins \{ n \in \partial \R_+^{J_t} \tq \cg \nu , n \cd \le 0 \} \p$$
Donc $p(V \cap \uni_{\alpha \in J_t} \alpha^\ort)$ est un espace contractile sauf s'il est de la forme :
$$ (\Diamond) \;\;  \partial \R_+^{J_t} \moins \{0\} \p$$
Dans ce cas, $p(V \cap \uni_{\alpha \in J_t} \alpha^\ort)$ est hom\'eomorphe \`a $\R^{|J_t|-1} \moins \{0\}$ et :
$$\ti{H}^{j-1} ( V \cap \uni_{\alpha \in J_t} \alpha^\ort ) = \left\{\begin{array}{ll}
0 & \si j-1 \not= |J_t| - 2 \\
\R & \si j-1 = |J_t| - 2
\end{array}
\right.$$
(\cf par exemple \cite[chap. 4, th\'eor\`eme 6]{Spa:at}).

On est dans cette situation $(\Diamond)$ si et seulement si :
$$\{ n \in \partial \R_+^{J_t} \tq \cg \nu , n \cd \le 0 \} = \{0\}$$
et on v\'erifie que c'est \'equivalent \`a : $\nu \in Q_t \p$

Finalement, pour tout $t \in W \moins \{1\}$, tout $i \ge 0$ et tous $\lambda ,\mu$, on retrouve que :
$$H^{i - 2l(t) - 1} ( V(\lambda , t * \mu) , V(\lambda , t * \mu) \cap \uni_{\alpha \in J_t} \alpha^\ort ) = \R$$ si $i = 2 l(t) + |J_t| \et t * \mu \in \lambda + Q_t$ et que  dans tous les autres cas :
$$H^{i - 2l(t) - 1} ( V(\lambda , t * \mu) , V(\lambda , t * \mu) \cap \uni_{\alpha \in J_t} \alpha^\ort ) = (0) \p$$ \qed

\begin{rems}[et exemples]\label{rem:magni}
\begin{itemize}
\item[1.---] Pour tout poids entier dominant $\mu$ et tout faisceau inversible $\call{L}$ sur $\adh{G_{ad}}$, la multiplicit\'e du $\GXGT-$module
$$\Plus_{i \ge 0} H^i(\adh{G_{ad}} , \call{L})$$
selon $\End(L(\mu))$ est major\'ee par l'ordre de $W$. 

\item[2.---] Quel que soit le poids entier dominant $\mu$, l'ensemble des  degr\'es $i$ en lesquels $H^i(\adh{G_{ad}},\call{L})$ a une multiplicit\'e non nulle selon $\End(L(\mu))$, pour au moins un faisceau inversible $\call{L}$ est exactement :
$$\{ 2l(t) + |J_t| \tq t \in W\} \p$$ 

On en d\'eduit, par exemple, que si $i = 1 , 2 ,$ ou $4$, alors :
$$ \; H^i(\adh{G_{ad}} , \call{L}_\lambda) = (0) $$
pour tout groupe $G$ et tout caract\`ere $\lambda$.

\item[3.---] Soit $\{ \omega_\alpha \tq \alpha \in \Delta \}$ la base des poids fondamentaux. Posons pour toute partie  $J$ de  $\Delta$ :
 
$$P_J := \left\{ \sum_{\alpha \in \Delta} p_\alpha \omega_\alpha \tq \qq \alpha ,\, p_\alpha \in \Z \et [ p_\alpha < -1 \equi \alpha \in J ] \right\}$$
$$ Q_J := \left\{ \sum_{\alpha \in \Delta} q_\alpha \alpha \tq \qq \alpha, \, q_\alpha \in \Z \et [ q_\alpha >0 \equi \alpha \in J] \right\} \p$$

Avec ces notations, si $\mu$ est dominant, la bijection $w \donne w * \mu$ de $W$ sur $W * \mu$ induit une bijection :
$$\left\{ t \in W \tq t *  \mu \in \lambda + Q_t  \right\} \sta{\sim}{\to} \uni_{J \sub \Delta} \left\{ \nu \in W * \mu \cap (\lambda + Q_J) \cap P_J \right\} \p$$
 
Si $\nu$ est un poids entier tel que $\nu + \rho$ est r\'egulier, alors on note $\nu^+$ l'unique poids dominant de $W * \nu$. On d\'eduit de la bijection ci-dessus que l'ensemble des poids entiers $\nu$  tels que  $\End(L(\nu^+))$ apparaisse dans la d\'ecomposition de la repr\'esentation :
$$\Plus_{i \ge 0} H^i(\adh{G_{ad}} , \call{L}_\lambda)$$
est exactement l'ensemble :
$$\uni_{J \sub \Delta} (\lambda + Q_J) \cap P_J \p$$

\item[4.---] \label{rem:magni4} On obtient aussi que, pour tous les  groupes adjoints $G$ de rang $2$, les multiplicit\'es $m^i_\lambda(\mu)$ sont $0$ ou $1$. 

En effet, soit $\Delta := \{\alpha, \beta\}$ ; soient $w_1,w_2 \in W$ v\'erifiant :
$$\left\{\begin{array}{c}
w_j * \mu \in \lambda +Q_{w_j} \\
2l(w_j) + |J_{w_j}| =i
 \end{array}
\right.
$$
pour $j =1 ,\, 2$ et un $0 < i < \dim G$ (si $i = 0$ ou $\dim G$, on sait d\'ej\`a que les multiplicit\'es $m^i_\lambda(\mu)$ sont $0$ ou $1$).

On a alors $J_{w_j} = \{\alpha\}$ ou $\{\beta\}$. Si $J_{w_1} = J_{w_2}$, alors $l(w_1) = l(w_2)$. On v\'erifie que cela entra\^{\i}ne : $w_1 = w_2$.

Si $J_{w_1} \not= J_{w_2}$, alors on a :
$$\lambda + Q_{\{\alpha\}} \cap P_{\{\alpha\}} \not= \vide \et \lambda + Q_{\{\beta\}} \cap P_{\{\beta\}} \not= \vide \p$$
Mais cela est impossible pour chacun des syst\`emes $A_1 \croi A_1, A_2, B_2, G_2$.
 
\item[5.---] On d\'eduit aussi de la formule du th\'eor\`eme \ref{thm:magni} que pour tout $G$ et  tout  caract\`ere $\lambda \in \ti{\cal X}$, les repr\'esentations $H^3(\adh{G_{ad}} , \call{L}_\lambda)$ sont sans multiplicit\'e. En effet, on a :
$$2l(w) + |J_w| = 3 \equi \exist \alpha \in \Delta , \; w=s_\alpha \p$$
De plus, si $\mu \in \ti{\cal X}^+$, $\lambda \in \ti{\cal X}$ et $\alpha, \beta \in \Delta$, alors :
$$s_\alpha * \mu \in \lambda + Q_{s_\alpha} \et s_\beta * \mu \in \lambda + Q_{s_\beta} \impliq s_\alpha * \mu - s_\beta * \mu \in Q_\alpha - Q_\beta \; (**)\p$$
Mais comme : $\ma s_\alpha * \mu - s_\beta * \mu =  \underbrace{-(\cg \mu ,\alpha^\ch \cd + 1 )}_{<0} \alpha + \underbrace{(\cg \mu ,\beta^\ch \cd + 1)}_{>0} \beta $ et puisque : $\ma Q_\alpha - Q_ \beta \sub \Z_{>0} \alpha - \Z_{>0} \beta + \Z (\Delta \moins \{ \alpha,\beta\})$, il faut que $\alpha = \beta$ pour que la condition $(**)$ soit v\'erifi\'ee.

\item[6.---] Les multiplicit\'es $m^i_\lambda(\mu)$  peuvent n\'eanmoins \^etre diff\'erentes de $0$ ou $1$ : si par exemple, $G = PSO(8, \C)$ (c'est le type $D_4$), on peut montrer \`a l'aide de 3.--- que pour $i=5$ :
$$ \qq \mu \in \ti{\cal X}^+ , \, \exist  \lambda \in \ti{\cal X} , \;  m^5_\lambda(\mu) = 3 \p$$

\item[7.---] Contrairement au cas du th\'eor\`eme de Borel-Weil-Bott, il peut arriver arriver que $m^i_\lambda(\mu) >0$ pour plus d'un degr\'e $i$ ($\lambda$ et $\mu$ \'etant fix\'es). Par exemple, lorsque $G$ est le groupe adjoint de type $F_4$, on obtient gr\^ace au th\'eor\`eme \ref{thm:magni} qu'il existe une infinit\'e de caract\`eres $\lambda \in \ti{\cal X}$ tels que $m_\lambda^{10}(0) > 0$ et $m^{11}_\lambda(0) > 0$. 

\end{itemize}
 \end{rems}
\tri

Maintenant, on va d\'emontrer le th\'eor\`eme \ref{thm:princi}. 

{\it On aura besoin de $\goth g$, l'alg\`ebre de Lie de $G$. On notera aussi $U(\goth g)$ son alg\`ebre enveloppante.}

On proc\`ede en trois grandes \'etapes. On commencera dans la section qui suit par rappeler la d\'efinition des groupes de cohomologie \`a support et le complexe de Grothendieck-Cousin ; ceux-l\`a apparaissant dans celui-ci. Dans le contexte du fibr\'e en droites $\call{L}_h$ sur la compactification $\G$, ce complexe a un double int\'er\^et. Non seulement son homologie est exactement la cohomologie $H^*(\G ,\call{L}_h)$ mais, ses termes, qui sont des groupes de cohomologie \`a support, sont aussi des repr\'esentations de $\goth g \croi \goth g$ dont on peut d\'ecrire pr\'ecis\'ement les sous-quotients simples de dimension finie, avec leurs multiplicit\'es. On va \'etudier ces groupes de cohomologie au cours de la deuxi\`eme \'etape. Pour cela, on ne perdra rien en se pla\c{c}ant dans le cadre plus g\'en\'eral des vari\'et\'es r\'eguli\`eres. Gr\^ace \`a cette \'etude, on pourra en la derni\`ere \'etape \'eliminer la plupart des termes du complexe de Grothendieck-Cousin. Pour ceux qui resteront, on se ram\`enera d'abord \`a la d\'etermination de groupes de cohomologie \`a support dans des vari\'et\'es qui sont $\ll$ presque des espaces affines $\gg$. Enfin pour conclure, on r\'eduira tout \`a un calcul dans le contexte torique. 
\section{Cohomologie \`a support}\label{sec:tGC}

Avec la d\'efinition des groupes de cohomologie \`a support, on redonne un r\'esultat d'annulation qu'on utilisera souvent.
\subsection{D\'efinition et th\'eor\`eme d'annulation}

Soient $X$ un espace topologique et $Z_2 \sub Z_1 \sub X$ deux ferm\'es de $X$. Soit ${\cal F}$ un faiseau ab\'elien sur $X$.

\begin{defi}
On note $\Gamma_{Z_1/Z_2}({\cal F})$ le quotient de groupes ab\'eliens :
$$\{\sigma \in {\cal F}(X) \tq \sigma \! \mid_{X \moins Z_1} = 0 \} \: / \: \{\sigma \in {\cal F}(X) \tq \sigma \! \mid_{X \moins Z_2} = 0\} \p$$
On appelle alors $\ll$ $i-$\`eme groupe de cohomologie de ${\cal F}$ \`a support dans $Z_1/Z_2$ $\gg$, et on  note $H^i_{Z_1/Z_2}({\cal F})$, le $i-$\`eme groupe d\'eriv\'e \`a droite du foncteur :
$${\cal F} \donne \Gamma_{Z_1/Z_2}({\cal F}) \p$$
Si $Z$ est seulement localement ferm\'e dans $X$, on d\'efinit : $\ma H^i_Z({\cal F}) :=  H^i_{\adh{Z} / \adh{Z} \moins Z}({\cal F}) $.
 \end{defi}

Cette notion g\'en\'eralise celle de cohomologie $\ll$ tout court $\gg$, en effet : $\ma H^i_X({\cal F}) = H^i(X, {\cal F})$. \`A l'oppos\'e, pour tout $i$ : $\ma H^i_\vide ({\cal F}) = H^i_{Z/Z}({\cal F}) = (0) $.

 Dans le contexte des vari\'et\'es on a :
\begin{thm}[{\cite[th. 9.5 et 9.6]{Ke:GCc}}]\label{lem:annul}
Soit $X$ une vari\'et\'e lisse et irr\'eductible, $Z$ une sous-vari\'et\'e affine de $X$, lisse et irr\'eductible. Alors, pour tout faisceau ${\cal F}$ coh\'erent et localement libre sur $X$, on a :
$$ \qq i \not= \codim(Z,X) , \; H^i_Z({\cal F}) = (0) \p$$
 \end{thm}

\begin{rem}
La condition $\ll$sous-vari\'et\'e affine, lisse et irr\'eductible$\gg$ est par exemple v\'erifi\'ee par les $B-$orbites de $X$, pour un groupe connexe et r\'esoluble $B$ qui op\`ere sur $X$.
 \end{rem}

\tri 

Rappelons aussi que pour une $G-$vari\'et\'e $X$, pour une sous-vari\'et\'e quelconque $Z$ de $X$ et pour un faisceau $\ti{G}-$lin\'earis\'e sur $X$, G. Kempf a montr\'e que les groupes de cohomologie $H^i_Z({\cal F})$ sont naturellement des $\goth g-$modules (\cf \cite[11.1, 11.3, 11.6]{Ke:GCc}). 

\subsection{Th\'eor\`eme de Grothendieck-Cousin}\label{ssec:tGC}

Ce th\'eor\`eme met en lumi\`ere le r\^ole de la cohomologie \`a support. 

\begin{thm}[\cite{MR:gr},\cite{Bo:gc},\cite{Ke:GCc},\cite{Ha:Rd}]\label{thm:cGC}
Soit $X$ un espace topologique. Soit $X \con Z_0 \con Z_1 \con \sus \con Z_n \con Z_{n+1} = \vide$ une filtration de $X$ par des sous-espaces ferm\'es. 
\begin{liste}
\item[a)] Pour tout faisceau ${\cal F}$ ab\'elien sur $X$, on a un complexe, le $\ll$ \re {complexe de Grothendieck-Cousin} $\gg$ :
$$0 \to  H^0_{Z_0/Z_1}({\cal F}) \sta{d^0}{\to} H^{1}_{Z_1/Z_2}({\cal F}) \sta{d^1}{\to} \sus \sta{d^{n-1}}{\to} H^{n}_{Z_n}({\cal F}) {\to} 0 \p$$
\item[b)] Si de plus la condition suivante est v\'erifi\'ee :
$$\qq p , \: \qq q \not= 0,\; H^{p+q}_{Z_p/Z_{p+1}} ( {\cal F} ) = (0) \; \; (\mbox{\ding{86}})\;\;,$$
alors pour tout $i \ge 0$, $H^i_{Z_0}(X, {\cal F})$ est le $i$-\`eme groupe d'homologie du complexe :
$$0 \to H^0_{Z_0/Z_1}({\cal F}) \sta{d^0}{\to} H^{1}_{Z_1/Z_2}({\cal F}) \sta{d^1}{\to} \sus \sta{d^{n-1}}{\to} H^{n}_{Z_n}({\cal F}) \sta{d^n}{\to} 0 \p$$
 \end{liste}

\end{thm}

\tri

D'apr\`es le th\'eor\`eme \ref{lem:annul}, la condition $(\mbox{\ding{86}})$ est v\'erifi\'ee quand les trois suivantes sont remplies :
\begin{liste}
\item $X$ est une vari\'et\'e alg\'ebrique lisse et irr\'eductible, 
\item ${\cal F}$ est un faisceau coh\'erent et localement libre sur $X$,
\item pour tout $p \ge 0$, la sous-vari\'et\'e $Z_p \moins Z_{p+1}$ est vide, ou est lisse et de codimension pure \infra{
\cad que toutes les composantes irr\'eductibles ont la m\^eme codimension.}
$= p$ dans $X$.
 \end{liste}


\tri

Appliquons ce th\'eor\`eme au faisceau inversible $\call{L}_h$ et \`a la compactification r\'eguli\`ere $\G$. On a donc besoin d'une filtration de $\G$ par des ferm\'es. Notons $Z_p$ la r\'eunion des $B \croi B^--$orbites de $\G$ de codimension $\ge p$ ($p$ est un entier positif). 
On consid\'erera la filtration :

$$\G = Z_0 \con Z_1 \sus \con  \sus $$
Puisque la vari\'et\'e $\G$ n'a qu'un nombre fini de $B \croi B^--$orbites (\cf  \cite[\S 2.1]{B:Inf}) et que l'adh\'erence d'une orbite est une union d'orbites de codimension plus grande, pour tout $p \ge 0$, les $Z_p$ sont des ferm\'es.
\tri

On va analyser le complexe de Grothendieck-Cousin qui appara\^{\i}t :

$$0 \to H^0_{Z_0/Z_1}(\call{L}_h) \to H^1_{Z_1/Z_2}(\call{L}_h) \to \sus$$

C'est un complexe de $\goth g \croi \goth g-$modules dont le $p-$i\`eme terme est :
$$H^p_{Z_p/Z_{p+1}} (\call{L}_h) = H^p_{Z_p - Z_{p+1}} (\call{L}_h) = \Plus_{ \Omega } H^p_{\Omega} (\call{L}_h) \;,$$
cette somme se faisant sur les $B \croi B^--$orbites $\Omega$ de $\G$ de codimension $p$. De plus, comme le faisceau $\call{L}_h$ est $\GXGT-$lin\'earis\'e sur $\G$, pour toutes les $B \croi B^-$orbites $\Omega$ de $\G$, l'action de l'alg\`ebre de Lie de $B \croi B^-$ s'int\`egre en une action (rationnelle) du groupe $\ti{B} \croi \ti{B^-}$. On dit que les groupes $H^p_\Omega(\call{L}_h)$ sont des $\goth g \croi \goth g-\ti{B} \croi \ti{B^-}-$modules (\cf \cite[pp. 373, 374, 384]{Ke:GCc}). Pour les \'etudier, on va en donner des filtrations. C'est l'objet de la partie qui suit.

\section{Filtration des groupes de cohomologie \`a support}\label{sec:figcs}

{\it Dans toute cette partie, $X$ sera une vari\'et\'e compl\`ete et r\'eguli\`ere (pour le groupe $G$).}
 
Soit $L$ un fibr\'e en droites $G-$lin\'earis\'e sur $X$. On va filtrer les groupes de cohomologie de $L$ \`a support dans une $B-$orbite $\Omega$ de $X$. Notre but est d'obtenir, comme gradu\'es associ\'es des $\goth g-$modules tr\`es particuliers : des groupes de cohomologie de fibr\'es en droites sur un espace homog\`ene $G/H$ et \`a support dans une $B-$orbite $BH/H$, avec un sous-groupe $H$ de $G$ contenant $T$. Pour ces groupes de cohomologie, le support $BH/H$ est un espace affine et on peut calculer leur caract\`ere comme $T-$modules. De plus, lorsque $G/H$ est une vari\'et\'e projective (\cad une vari\'et\'e de drapeaux), ces modules sont des repr\'esentations famili\`eres de $\goth g$ : des modules de Verma g\'en\'eralis\'es. 

En fait, on ne s'attaque pas directement aux groupes de cohomologie \`a support dans les $B-$orbites. On \'etudiera d'abord la cohomologie \`a support dans les cellules de Bialynicki-Birula (dont on rappelle la d\'efinition ci-dessous) : pour obtenir les filtrations voulues, c'est plus commode.  Ensuite, on passera au cas o\`u le support est une $B-$orbite. Enfin, on verra qu'on aboutit au cas \'evoqu\'e ci-dessus, avec $G/H$ projective, lorsque tous les $T-$points fixes de $X$ sont dans une $G-$orbite ferm\'ee (comme on l'a vu, cette derni\`ere condition est satisfaite par la compactification r\'eguli\`ere $\G$).

\subsection{Les cellules et les orbites des vari\'et\'es r\'eguli\`eres}\label{sec:vr}

Dans les vari\'et\'es r\'eguli\`eres, les cellules de Bialynicki-Birula jouent un r\^ole remarquable ; elles permettent notamment de param\'etrer les $B-$orbites et, on les utilisera aussi pour les filtrations de groupes de cohomologie \`a support.

\tri

L'ensemble des points fixes de $T$ dans $X$ est fini.  Notons le $X^T$. Il existe un sous-groupe \`a un param\`etre $\zeta$ (de $T$), dominant (\cad tel que pour tout $\alpha \in \Delta$, $\cg \alpha ,\zeta \cd \ge 0$) et dont les points fixes, dans $X$, sont aussi fix\'es par $T$ (\cf par exemple \cite[lemme 2.3]{BB:aag}). 

{\it Fixons pour toute cette partie un tel $\zeta$.}

 Pour chaque $x \in X^T$, la cellule de Bialynicki-Birula 
$$X(x) := \{ y \in X \tq \limz \zeta(a) . y =x\}$$
est une sous-vari\'et\'e de $X$, isomorphe \`a un sous$-T-$module de $T_xX$, l'espace tangent \`a $X$ en $x$ (\cf \cite {BB:aag}). Si $Z$ est une sous-vari\'et\'e ferm\'ee et $G-$invariante de $X$, on notera $Z(x) := X(x) \cap Z$ la cellule de Bialynicki-Birula de $Z$ associ\'ee \`a $x$ et \`a $\zeta$. Ces sous-vari\'et\'es $X(x)$ (et $Z(x)$) sont $B-$invariantes car, pour tout $b \in B$, la limite $$\limz \zeta(a)b\zeta(a\inv)$$ existe et appartient \`a $T$. D\'esormais, on appellera simplement {\it cellules} les cellules de Bialynicki-Birula associ\'ees aux points de $X^T$ et \`a $\zeta$. 

\tri

On peut retrouver les $B-$orbites de $X$ gr\^ace aux cellules. En effet, d'apr\`es \cite[\S 2.1 p. 219 et pro. du \S 2.3]{BL:loc}, d'une part l'intersection d'une cellule et d'une $G-$orbite de $X$ est soit vide soit une $B-$orbite ; d'autre part, on obtient ainsi toutes les $B-$orbites de $X$.

\subsection{Filtration de la cohomologie \`a support dans les cellules}
{\it Soit ${\cal D}$ l'ensemble des diviseurs limitrophes de $X$} (ce
sont les composantes irr\'eductibles de $X \moins X^0_G$, \cf la
d\'efinition \ref{def:vreg}, p. \pageref{def:vreg}).
Soit $x \in X^T$. Pour \'etudier les groupes de cohomologie \`a support dans les cellules, notre m\'ethode consiste \`a filtrer selon l'ordre d'annulation le long des diviseurs limitrophes. Et, pour se ramener \`a de la cohomologie sur $G \hp x$ et \`a support dans $B \hp x$, on utilisera notamment que : $X (x) \cap \adh{G \hp x} = B \hp x$. 

Soit $Z$ une sous-vari\'et\'e ferm\'ee, $G-$invariante et irr\'eductible, de $\adh{G.X(x)}$ qui contient $x$. On notera $z(x)$ la codimension de la cellule $Z(x)$ dans $X$ et $b(x)$ celle de l'orbite $B \hp x$ dans l'orbite $G \hp
x$. D'apr\`es le th\'eor\`eme \ref{lem:annul}, pour tout faisceau ${\cal M}$ localement libre et de rang fini  sur $X$, on a : $\ma
H^i_{Z(x)}({\cal M}) = (0)$ si $i \not= z(x)$. Le th\'eor\`eme suivant
concerne le groupe de cohomologie \`a support $\ma H^{z(x)}_{Z(x)}({\cal M})$.
 
\begin{thm}\label{thm:filtcel}
 Soit ${\cal M}$ un faisceau coh\'erent, localement libre et $\ti{G}-$lin\'earis\'e sur $X$. Pour tous $m \ge -1 ,\,n \ge 0$, on note ${\cal M}^m_n$ le faisceau $\ma \Plus_{E} {\cal M}(E)$ o\`u $E$ d\'ecrit l'ensemble des diviseurs de la forme :
$$ \sum_{D \in {\cal D} \atop Z  \sub D} (m_D + 1) D - \sum_{ D \in {\cal D} \atop x \in D , Z \not \sub D} n_D D$$
avec des entiers $m_D,n_D \ge 0$ v\'erifiant :
$\sum_D m_D = m \et  \sum_D n_D =n$.

Alors, avec ces notations, on a une filtration de $\goth g-$modules : $$\ma H^{z(x)}_{Z(x)}({\cal M}) = \uni_{m \ge -1 ,\, n \ge 0} F^m_n \;,$$ d\'ecroissante selon l'indice $n$, telle que, pour tous $m,n \ge 0$, on ait : $\ma F^{-1}_n = (0) \,,\; \inter_{n \ge 0} F^{m}_n = F^{m-1}_0$ et dont les quotients successifs sont :
$$F^m_n / F^m_{n+1} = H^{b(x)}_{B \: \cdot \: x}({\cal M}^m_n\res{G \: \cdot \:x}) \p$$

 \end{thm}

\begin{rems}
\begin{itemize}
\item[1.---] Dans la suite,  on appliquera ce th\'eor\`eme \`a des
faisceaux inversibles sur $X$.
\item[2.---] Si on supprime l'hypoth\`ese $\ll$ $\ti{G}-$lin\'earis\'e
$\gg$, on perd seulement le caract\`ere $\goth g-$\'equivariant de la
filtration.
\item[3.---] On a not\'e pour tout diviseur $D$ de $X$, ${\cal M}(D)$
le faisceau $\ma {\cal M}\tenso{\fstr_X}\fstr_X(D)$.
 \end{itemize}
 \end{rems}

\begin{dem}
Pour tout sch\'ema $X'$, le signe ${}^\ch$ symbolisera le dual des $\fstr_{X'}-$modules. On emploiera aussi les notations suivantes : pour tous ferm\'es $Z_1 \sub Z_2$ de $X$,  \begin{itemize}
\item ${\cal I}_{Z_1/Z_2}$ le faisceau d'id\'eaux de d\'efinition de $Z_1$ dans $Z_2$ (et parfois, pour abr\'eger : ${\cal I}_{Z_1} : = {\cal I}_{Z_1/X}$) ;
\item ${\cal N}_{Z_1/Z_2} := \left({\cal I}_{Z_1/Z_2} / {\cal I}^2_{Z_1/Z_2} \right)^\ch$ le faisceau normal de $Z_1$ par rapport \`a $Z_2$ (\cf \cite[d\'ef. du \S 8.19]{Ha:AG}) ;
\item  $\omega_{Z_1/Z_2} := \Vect^{z_2 - z_1} {\cal N}_{Z_1 / Z_2}$ le faisceau canonique de $Z_1$ par rapport \`a $Z_2$, lorsque $Z_1$ et $Z_2$ sont lisses et irr\'eductibles, de dimensions respectives $z_1$ et $z_2$  (\cf \cite[d\'ef. b) p. 141]{Ha:Rd}).
\end{itemize}

On notera enfin $S^p$ la puissance sym\'etrique $p-$i\`eme. 
\tri 
Remarquons pour commencer que, d'apr\`es \cite [th. du \S 1.4, ii)]{B:Inf}, $\adh{G.Z(x)} = Z$ et que, d'apr\`es \cite[d\'ebut du \S 2.4]{BB:alrrv}, on a des d\'ecompositions en faisceaux inversibles :
$${\cal N}_{Z/X} \res{ \adh{G \hp x}} = \Plus_{D \in {\cal D} , D \con Z} \fstr_X(D) \res{ \adh{G \hp x}} $$
$$ {\cal N}_{ G \hp x / Z}\res{\adh{G \hp x}} = \Plus_{ D \in {\cal D} \atop D \ni x , D \not \con Z} \fstr_X(D)\res{\adh{G \hp x}}$$
$$\omega_{Z/X}\res{\adh{G \hp x}} = \fstr_X \left(  \sum_{D \in {\cal D}, D \con Z} D \right) \res{\adh{G \hp x}} \p$$

\noindent 1.--- {\bf Cas o\`u $Z = X$}

On a alors : $\ma {\cal M}^0_n \res{\adh{G \hp x}} = {\cal M} \tens S^n{\cal N}_{\adh{G \hp x} / X}^\ch$. On va filtrer selon l'ordre d'annulation le long de $\adh{G \hp x}$ : on remarque, effectivement, que : $\ma {\cal M} = \uni_{n \ge 0} {\cal M} \tens {\cal I}_{\adh{G \hp x}}^n$ d'o\`u :

$$ H^{z(x)}_{Z(x)} ({\cal M}) = \limd{n \ge 0} H^{z(x)}_{Z(x)}(  {\cal M} \tenso{ {\fstr}_X } {\cal I}_{\adh{G \hp x}}^{n} ) \p$$
On va voir que cette limite directe est en fait une r\'eunion (croissante) et on va en d\'eterminer les quotients successifs. Pour tout $n \ge 0$, on a une suite exacte courte :

$$
0 \to {\cal M} \tenso{ {\fstr}_X } {\cal I}_{\adh{G \hp x}}^{n+1} \to {\cal M} \tenso{ {\fstr}_X } {\cal I}_{\adh{G \hp x}}^{n} \to {\cal M} \tenso{ {\fstr}_X } S^n {\cal N}_{\adh{G \hp x}/X}^\ch \to 0
$$
d'o\`u l'on d\'erive une suite exacte :
$$(*) \;\; 0 \to H^{z(x)}_{Z(x)}( {\cal M} \tens {\cal I}_{\adh{G \hp x}}^{n+1} ) \to$$
$$  H^{z(x)}_{Z(x)}(  {\cal M} \tens {\cal I}_{\adh{G \hp x}}^{n} ) \to  H^{b(x)}_{ B\hp x}({\cal M} \tens S^n {\cal N}_{\adh{G \hp x} / X}^\ch ) \to 0 \p$$

En effet, comme $Z(x)$ est une sous-vari\'et\'e affine et lisse de $X$, d'une part, on d\'eduit du th\'eor\`eme \ref{lem:annul} que pour $i \not = z(x)$, $\ma H^i_{Z(x)} ({\cal M} \tens {\cal I}_{\adh{ G \hp x}}^n ) = (0)$ et, d'autre part, comme $Z(x)$ et $\adh{G \hp x}$ se coupent proprement dans $Z = X$ ( \cf \cite[ii) du th. du \S 1.4]{B:Inf}), avec $Z(x) \cap \adh{G \hp x} = B \hp x$, et ${\cal N}_{\adh{G \hp x}/X}$ est localement libre sur $\adh{G \hp x}$, on a de m\^eme : 
$$\qq i \not= z(x), \; H^{i}_{ Z(x)}({\cal M} \tens S^n {\cal N}_{\adh{G \hp x} / X}^\ch ) = H^{i}_{ B \hp x}({\cal M} \tens S^n {\cal N}_{\adh{G \hp x} / X}^\ch ) = (0) \p$$

Pour terminer, en se servant de $(*)$, il ne reste plus qu'\`a poser : $\ma F^0_n := H^{z(x)}_{Z(x)}({\cal M} \tens {\cal I}_{\adh{G \hp x}}^n)$ pour $n \ge 0$. 
On v\'erifie aussi que : 
$$\inter_{ n \ge 0} H^{z(x)}_{Z(x)} ({\cal M} \tens {\cal I}_{\adh{G \hp x}}^n ) = (0) \p$$
En effet, on peut majorer les caract\`eres des groupes de cohomologie \`a support dans une cellule. On montre ainsi que les poids du tore $\ti{T}$ dans le $\ti{T}-$module $H^{z(x)}_{Z(x)} ({\cal M} \tens {\cal I}_{\adh{G \hp x}}^n )$ se d\'ecalent quand $n$ cro\^{\i}t et que chacun finit par dispara\^{\i}tre quand  $n$ tend vers l'infini (\cf \cite[th. II.3.2]{these}). 

\tri

\noindent 2.--- {\bf On se ram\`ene au cas de l'\'etape pr\'ec\'edente}

D'apr\`es \cite[pro.2.5]{BB:alrrv}, la vari\'et\'e $Z$ (comme toute adh\'erence de $G-$orbite) est encore r\'eguli\`ere. Notons $z$ sa codimension dans $X$.

 Puisque $Z(x) \sub Z$, gr\^ace \`a \cite[th. 2.8]{Gr:LC} et \`a \cite[lem. 8.5.d)]{Ke:GCc}, on a :
$$H^{z(x)}_{Z(x)}({\cal M}) = \limd{m \ge 0} H^{z(x) - z}_{Z(x)} \left( \Ext^{z}_{\fstr_X} ( \fstr_X / {\cal I}^m_{Z} , {\cal M} ) \right) \p$$
Cette limite directe est encore une r\'eunion croissante. En effet, on tire de la suite exacte courte :
$$ 0 \to S^m {\cal N}^\ch_{Z /X} \to \fstr_X / {\cal I}^{m+1}_{Z} \to \fstr_X/ {\cal I}_{Z}^m \to 0$$
une autre suite exacte de faisceaux :
 $$
0 \to \Ext^{z}_{\fstr_X} ( \fstr_X / {\cal I}^m_{Z} , {\cal M} ) \to \Ext^{z}_{\fstr_X} ( \fstr_X / {\cal I}^{m+1}_{Z} , {\cal M} ) \to \Ext^{z}_{\fstr_X} (  S^m{\cal N}_{Z/X}^\ch , {\cal M} ) \to 0 \p$$
Or, puisque le faisceau $S^m{\cal N}_{Z/X}^\ch$ est localement libre sur $Z$, on a, d'apr\`es \cite[pro. III.7.2]{Ha:Rd} :
$$\Ext^{z}_{\fstr_X} (  S^m{\cal N}_{Z/X}^\ch , {\cal M} ) = \omega_{Z/X} \tens S^m{\cal N}_{Z/X}\tens {\cal M} $$
qui est un faisceau localement libre sur $Z$.

Donc, en posant $\ma F^m := H^{z(x) - z}_{Z(x)} \left( \Ext^{z}_{\fstr_X} (  \fstr_X / {\cal I}^m_{Z} , {\cal M} ) \right)$, on obtient une filtration croissante $ \ma H^{z(x)}_{Z(x)}({\cal M}) = \uni_{m \ge 0} F^m$ de quotients succesifs : 
$$F^{m+1} / F^{m} = H^{z(x) - z} _{Z(x)}(\omega_{Z/X} \tens S^m{\cal N}_{Z/X} \tens {\cal M}) \p$$
On peut alors appliquer la 1\`ere \'etape \`a la vari\'et\'e r\'eguli\`ere $Z$ et au faisceau localement libre : $\ma \omega_{Z/X} \tens S^m{\cal N}_{Z/X} \tens {\cal M}$. On trouve pour tout $m \ge -1$, une filtration d\'ecroissante $ F^{m+1}/F^{m} = \uni_{n \ge 0} \adh{F}^m_n$. Puisque ${\cal M}^m_n\res{\adh{G \hp x}} =  \omega_{Z/X} \tens S^m{\cal N}_{Z/X} \tens {\cal M} \tens S^n {\cal N}_{\adh{G \hp x} / Z}^\ch$, on a :
$$\qq m ,n,  \; \adh{F}^m_n / \adh{F}^m_{n+1} = H^{b(x)}_{ B \hp x} ({\cal M}^m_n \res{ G \hp x}) \p$$
On conclut la d\'emonstration en prenant pour $F^m_n$ le relev\'e de $\adh{F}^m_n$ dans $F^{m+1}$ $(m \ge -1)$. 

\end{dem}

\tri

Soit $\Omega$ une $B-$orbite de $X$, et $b$ sa codimension. On va maintenant exprimer les groupes de cohomologie \`a support dans $\Omega$ \`a l'aide de groupes de cohomologie \`a support dans une certaine cellule. Pour cela, on va utiliser que l'orbite $\Omega$ est l'intersection de la $G-$orbite $G . \Omega$ et d'une cellule $X(x)$. Ici, le point $x$ est la limite $\ma \lim_{a \to 0} \zeta(a) .y$ pour un $y$ quelconque de $\Omega$ (cette limite est ind\'ependante du $y$ choisi).

Notons $Z_\Omega$ l'adh\'erence de la $G-$orbite $G . \Omega$ et $D$ la somme des diviseurs limitrophes contenant $x$ mais non $\Omega$. Dans le th\'eor\`eme pr\'ec\'edent, on a \'etudi\'e les groupes $H^b_{Z_\Omega(x)} ({\cal M})$ pour un faisceau localement libre ${\cal M}$ sur $X$. Ils apparaissent dans la filtration de la proposition suivante :


\begin{pro}\label{pro:filtrat}
Soit ${\cal M}$ un faisceau coh\'erent, localement libre et $\ti{G}-$lin\'earis\'e sur $X$.

Alors, on a une filtration croissante de $\goth g-$modules :

$$ H^b_\Omega({\cal M}) = \Uni_{ k \ge 0} H^b_{Z_\Omega(x)} ({\cal M}(kD)) \p$$
\end{pro}
\begin{dem}

Puisque $G.\Omega = Z_\Omega \moins D$, on a : $G.\Omega = Z_\Omega (x) \cap (X \moins D)$. Or si on note $j : X \moins D \inc X$ l'inclusion canonique, on a :
$$\limd{k \ge 0} {\cal M }(kD) = j_* {\cal M}\res{X \moins D} \p$$

Or, $j$ est un morphisme affine (il suffit de le v\'erifier pour chaque diviseur limitrophe, qui est lisse) ; d'o\`u : 
$$H^b_{\Omega}({\cal M}) = H^b_{j\inv(Z_\Omega(x))} ({\cal M}\res{X \moins D}) =  H^b_{Z_\Omega(x)}(j_* {\cal M}\res{X \moins D} ) = \limd{k \ge 0} H^b_{Z_\Omega(x)}({\cal M}(kD) )\p$$

En outre, les morphismes :
  
$$H^b_{Z_\Omega(x)}({\cal M}(kD) ) \to H^b_{Z_\Omega(x)}({\cal M}((k+1)D) )$$
sont injectifs pour tout $k \ge 0$ (car chaque composante irr\'eductible de $D$ coupe proprement $Z_\Omega(x)$ dans $Z_\Omega$).

Par cons\'equent, on a la filtration voulue.
\end{dem}

\subsection{Suites de composition}\label{sec:sdc}

Soit $\call{L}$ un faisceau $\ti{G}-$lin\'earis\'e sur $X$. On suppose encore que $\Omega$ est une $B-$orbite de codimension $b$ dans $X$.
En g\'en\'eral, les $\goth g-$modules $H^b_\Omega (\call{L})$ ne sont pas de type fini. N\'eanmoins, ils admettent parfois une d\'ecomposition naturelle en somme directe de $\goth g-$modules de longueur finie. C'est ce que nous allons maintenant \'etudier.

Soit $Z(\goth g)'$ l'ensemble des caract\`eres centraux du centre $Z(\goth g)$ de $U(\goth g)$. Rappelons (\cf \cite[\S 7.8.15]{Dix:AE}) que si $M$ est un $\goth g-$module et si $\chi \in Z(\goth g)'$, on peut d\'efinir, par r\'ecurrence, une suite croissante de sous-$\goth g-$modules de $M$ en posant :
$$M^0_\chi := (0) \, , \; \qq n >0 , \; M^n_\chi := \{ m \in M \tq \qq \z \in Z(\goth g) , \, \z.m - \chi(\z)m \in M^{n-1}_\chi \} \p$$

On appelle la r\'eunion croissante $\ma M_\chi : = \uni_{n \ge 0} M^n_\chi$ {\it l'espace propre g\'en\'eralis\'e de poids $\chi$ de $M$}. Remarquons que les $M_\chi$ sont en somme directe lorsque $\chi$ d\'ecrit $Z(\goth g)'$.

On reprend ici les notations du paragraphe \ref{sec:fi} sur $\ti{G} , \ti{B} , \etc .$ On supposera dans cette partie que $\call{L}$ est un faisceau $\ti{G}-$lin\'earis\'e. 

 Pour chaque caract\`ere central $\chi$ consid\'erons l'espace propre g\'en\'eralis\'e $\left( H^b_\Omega(\call{L}) \right)_\chi$.

Suivant \cite[pp. 373, 374, 384]{Ke:GCc}, les groupes $H^b_\Omega(\call{L})$ sont des $\goth g-\ti{B}-$modules ; \cad  que l'action de l'alg\`ebre de Lie de $B$ (par retriction de celle de $G$) $\ll$ s'int\`egre $\gg$ en une action rationnelle de $\ti{B}$. D'apr\`es \cite[\S 7.8.15]{Dix:AE} (\cf aussi \cite[annexe C1]{these}), il s'ensuit :

\begin{lem}
Comme $\goth g-$modules :
$$H^b_\Omega(\call{L}) = \Plus_{\chi \in Z(\goth g)'} (H^b_{\Omega}(\call{L}))_\chi \;\;\Box$$
 \end{lem}

On va montrer (proposition \ref{pro:sdc}) que les espaces propres g\'en\'eralis\'es $\left( H^b_\Omega(\call{L}) \right)_\chi$ sont des $\goth g-$modules de longueur finie lorsque le point $x$ fix\'e par $T$ et associ\'e \`a la $B-$orbite $\Omega$ (\cad l'unique $x$ tel que $\Omega = G . \Omega \cap X(x)$) appartient \`a une $G-$orbite ferm\'ee de $X$. Ce genre d'orbites $\Omega$ pr\'esente un avantage : dans les filtrations des groupes $H^b_\Omega (\call{L})$ obtenues \`a la section pr\'ec\'edente, ce sont des groupes de cohomologie de faisceaux sur l'orbite $G \hp x$ qui interviennent comme quotients successifs, et cette orbite est une vari\'et\'e de drapeaux. Ces groupes de cohomologie sont, dans ce cas, des $\goth g-$modules qui ont un caract\`ere central connu, c'est ce que l'on rappelle, entre autres, dans la section qui suit. Par cons\'equent, en fixant un caract\`ere central $\chi$, on obtiendra, gr\^ace au th\'eor\`eme et \`a la proposition \ref{thm:filtcel} et \ref{pro:filtrat}, une filtration explicite de l'espace propre g\'en\'eralis\'e correspondant : $\ma \left( H^b_\Omega(\call{L}) \right)_\chi$. Cette filtration se trouvera \^etre une suite de composition finie (\cf la proposition \ref{pro:sdc}).

Mais d'abord rappelons quelques $\goth g-\ti{B}-$modules particuliers :

\subsubsection{Les modules de Verma g\'en\'eralis\'es}

Soit $Q$ un sous-groupe parabolique de $G$ contenant $B^-$. Si $L \con T$ est le sous-groupe de Levi correspondant, de racines postives $\Phi^+_L$, on notera $$W^Q := \{ w \in W \tq \qq \alpha \in \Phi^+_L ,\: w(\alpha) >0 \}\p$$ Ainsi, les $B-$orbites de la vari\'et\'e de drapeaux $G/Q$ sont les :
$$BwQ/Q \tq w \in W^Q \p$$

\begin{rem}
Pour tout $w \in W^Q$, $l(w)$ est la codimension de $BwQ/Q$ dans $G/Q$.
 \end{rem}  

D'un autre c\^ot\'e, tout caract\`ere $\lambda$ de $\ti{Q}$ d\'etermine un unique faisceau $\call{L}$, inversible et $\ti{G}-$lin\'earis\'e sur la vari\'et\'e $G/Q$ (c'est celui pour lequel $\ti{Q}$ op\`ere via $\lambda$ dans la fibre $\call{L}\res{Q/Q}$). Nous le noterons $\call{L}_{G/Q}(\lambda)$. Le lemme suivant d\'ecrit les groupes de cohomologie \`a support $\ma H^{l(w)}_{BwQ/Q} (\call{L}_{G/Q}(\lambda))$ comme $\goth g-$modules. 

Pour l'\'enoncer, notons, pour tout caract\`ere $\lambda \in \ti{\cal X}$, $\chi_\lambda$ le caract\`ere central avec lequel le centre $Z(\goth g)$ op\`ere dans tous les $\goth g-$modules engendr\'es par un vecteur de plus haut poids $\lambda$ (\cf \cite[pro. 7.4.4]{Dix:AE}).  

\begin{thm}[{\cite[\S 3, pro. 3]{Bry:OD} et \cite[pro. 3.5]{BoBry:DO}}]\label{thm:BoBry}
Pour tout caract\`ere $\lambda$ de $\ti{Q}$, le $\goth g-\ti{B}-$module $$H^{l(w)}_{BwQ/Q} (\call{L}_{G/Q}(\lambda))$$
\begin{itemize}
\item admet $\chi_\lambda$ comme caract\`ere central ;
\item a une suite de Jordan-H\"{o}lder (finie) ;
\item a $w * \lambda$ comme plus haut poids.
 \end{itemize}
 \end{thm}

\begin{defi}
On notera $M^w_Q(\lambda)$ le $\goth g-$module $H^{l(w)}_{BwQ/Q}(\call{L}_{G/Q}(\lambda))$.
 \end{defi}

\begin{rem}
Il r\'esulte du th\'eor\`eme \ref{thm:BoBry} que le $\goth g-$module $M^w_Q(\lambda)$ a un $\ti{G}-$module comme sous-quotient simple si et seulement si $w * \lambda$ est dominant. De plus dans ce cas, le seul sous-quotient simple qui est un $\ti{G}-$module est $L(w* \lambda)$ et sa multiplicit\'e est $1$.
 \end{rem}

Si $Q = B^-$, $M^w_{B^-}(\lambda)$ est le module de Verma $w-$tordu $M^w(\lambda)$ de plus haut poids $w*\lambda$ d\'efini dans \cite[\S 2.2]{FF:AKM}.

\subsubsection{Quand le support est une orbite de rang maximal}\label{sec:rm}

La notion d'orbites de rang maximal est d\'efinie par exemple dans \cite[\S 3]{B:ooc}. En fait, ce sont exactement les $B-$orbites $\Omega$ de $X$ telles que la $G-$orbite $G \hp x$ du point limite :
$$x  = \lim_{a \to 0}\zeta(a). y \;\; (\qq y \in \Omega)$$
est ferm\'ee dans $X$ (\cf \cite[\S 3, th. 3]{B:ooc}). 

\begin{center}
\it 
On suppose dans ce paragraphe que $\Omega$ est une $B-$orbite de rang maximal.
 \end{center}

Pr\'ecisons quelques notations suppl\'ementaires avant d'\'enoncer et de d\'emontrer notre r\'esultat sur les suites de composition des $\goth g-$modules $\left( H^b_\Omega(\call{L}) \right)_\chi$. 

Soit $Q$ le sous-groupe parabolique de $G$ contenant $B^-$ associ\'e \`a $X$ : c'est le sous-groupe parabolique oppos\'e \`a $P$, le stabilisateur de la $B-$orbite ouverte de $X$ ; de plus, toutes les $G-$orbites ferm\'ees de $X$ sont isomorphes \`a $G/Q$ (\cf \cite[\S 2.2]{BB:alrrv} et \cite[\S 1.4]{B:Inf}).
 Puisque le point $x := \lim_{a \to 0} \zeta(a) . y$ (pour tout $y \in \Omega$) est dans une $G-$orbite ferm\'ee de $X$, il existe un unique $w \in W^Q$ tel que $w.x$ soit un point fixe de $Q$ (dans $X$). Notons $w_\Omega$ ce $w$ et $z_\Omega$ ce point $w_\Omega.x$. Les faisceaux inversibles sur $X$ sont d\'etermin\'es par leur restriction aux orbites ferm\'ees de $X$ et donc par leurs fibres en les points fixes de $Q$ dans $X$. Pour tout faisceau inversible et $\ti{G}-$lin\'earis\'e sur $X$, $\call{L}$, soit $p_\Omega(\call{L})$ le caract\`ere avec lequel le groupe $\ti{Q}$ op\`ere dans la fibre $\call{L}\res{z_\Omega}$.  

Avec ces notations, il est imm\'ediat que :
\begin{lem}\label{lem:dpdp}
D'une part :
$$b(x) := \codim ( B \hp x , G \hp x) = \codim (Bw_\Omega Q/Q , G/Q) =l(w_\Omega)$$
et d'autre part, pour tout faisceau $\call{L}$ inversible et $\ti{G}-$lin\'earis\'e sur $X$, on a l'isomorphisme de $\goth g-$modules :
 
$$ H^{b(x)}_{{B \hp x}} (G \hp x , \call{L}\res{G \hp x}) \iso M^{w_\Omega}_Q(p_\Omega(\call{L})) \;\;\Box$$
 \end{lem}

Pour un diviseur $B-$invariant $D$ de $X$, on notera $p_\Omega(D)$ le caract\`ere $p_\Omega(\fstr_X(D))$. 

\tri

Gardons les notations pr\'ec\'edentes. Si $\chi$ est un caract\`ere central, l'ensemble des caract\`eres $\nu \in \ti{\cal X}$ tels que $\chi_\nu = \chi$ est une classe de $\ti{\cal X}$ modulo l'action $*$ de $W$ (\cf \cite[pro. 7.4.7]{Dix:AE} ; notons la $W * \chi$. Du th\'eor\`eme et de la proposition sur les filtrations \ref{thm:filtcel} et \ref{pro:filtrat}, on d\'eduit la :
\begin{pro}\label{pro:sdc}
Soit encore ${\cal D}$ l'ensemble des diviseurs limitrophes de $X$. On consid\`ere les parties suivantes de l'ensemble des caract\`eres de $Q$ :
 $$
I_{\Omega} := \{ p_\Omega(D) \tq D \in {\cal D} , D \not\con \Omega \}$$
$$J_{\Omega} := \{ p_\Omega(D) \tq D \in {\cal D} , D \con \Omega\} \p$$

Alors, pour tout faisceau $\call{L}$, inversible et $\ti{G}-$lin\'earis\'e sur $X$ et pour tout caract\`ere central $\chi \in Z(\goth g)'$, le $\goth g-$module $\left( H^b_\Omega(\call{L}) \right)_\chi$ admet une suite de composition finie :

$$ \left(H^b_\Omega(\call{L})\right)_\chi = F_0 \supset F_1 \supset \sus \supset F_N \supset (0)$$
dont les quotients successifs sont \`a permutation pr\`es les :
$$M^{w_\Omega}_Q (\lambda) \tq \lambda \in \left( p_\Omega(\call{L}) + \Z I_{\Omega} + \Z_{ > 0} J_{\Omega} \right) \cap W * \chi \p$$
 \end{pro}

\begin{rem}
En particulier, les $\goth g-$modules $\ma  \left(H^b_\Omega(\call{L})\right)_\chi$ sont de longueur finie.
 \end{rem}

\begin{dem}

Gr\^ace au th\'eor\`eme \ref{thm:filtcel} et \`a la proposition \ref{pro:filtrat}, et avec leurs notations, on obtient une filtration infinie de $\goth g-$modules :
$$H^b_\Omega(\call{L}) = \uni _{m \ge -1, n,k \ge 0} F^{m,k}_n$$
dont les quotients successifs sont :
$$H^{b(x)}_{B \hp x} ({\cal M}^m_n(kD) \res{G \hp x}) \;\; (m \ge -1 , n,k \ge 0) \p$$

Lorsque $(m,n,k)$ varie dans $\Z_{\ge -1} \croi \Z_{\ge 0} \croi \Z_{\ge 0}$, \`a l'aide du lemme \ref{lem:dpdp}, on obtient comme quotients successifs les $\goth g-$modules :
$$
M^{w_\Omega}_Q(\lambda)$$
o\`u $\lambda$ d\'ecrit l'ensemble $p_\Omega(\call{L}) + \Z I_\Omega + \Z_{>0} J_\Omega$. (Remarquons que si $z_\Omega \notin D$, alors $p_\Omega(D) = 0$.) 

Maintenant, prendre la $\chi-$i\`eme composante $\ma \left( H^b_\Omega(\call{L}) \right)_\chi$ revient \`a ne garder, parmi ces $\goth g-$modules que ceux qui sont de caract\`ere central $\chi$ : \cad tels que : $\chi_\lambda = \chi$, ou encore $\lambda \in W * \chi$.

Le fait que chaque module de Verma g\'en\'eralis\'e appara\^{\i}t au plus une fois dans la suite de composition est une cons\'equence de l'ind\'ependance $\Z-$lin\'eaire des caract\`eres $p_\Omega(D)$, $D \in {\cal D}$ (\cf \cite[pro. A1, iv)]{B:Inf}). On obtient une suite de composition finie car chaque ensemble $W * \chi$ est fini.
\end{dem}

\tri

\subsubsection{Application}

En cons\'equence de la derni\`ere proposition, on trouve le r\'esultat suivant sur la structure des $\goth g \croi \goth g-$modules $H^{l(w) + l(w_0)} _{B w B^-} (\fstr_G)$ pour tout $w \in W$ :

\begin{pro}\label{pro:ep}
Soient $G$ est un groupe r\'eductif et connexe sur $\C$ (avec les notations de la p. \pageref{par:ici}) et $w$ un \'el\'ement du groupe de Weyl. Alors $H^i_{BwB^-}(\fstr_G) = (0)$, pour $i \not= l(w)+l(w_0)$, et les espaces propres g\'en\'eralis\'es non nuls de $H^{l(w) + l(w_0)} _{B w B^-} (\fstr_G)$ sont associ\'es aux caract\`eres centraux $\chi_{\lambda,-\lambda}$, $\lambda \in \ti{\cal X}$. Chaque tel espace propre $H^w_\lambda$ (associ\'e \`a $\chi_{\lambda,-\lambda}$) a une suite de composition finie de  $\goth g \croi \goth g-$modules :
$$H^w_\lambda = F_0 \supset F_1 \supset \sus \supset F_N \supset (0)$$
o\`u, \`a permutation pr\`es, les $F_i/F_{i+1}$ sont les :
$$M^w(\mu) \cten M( - \mu)^* \: , \; \mu \in W * \lambda \p$$
\end{pro}

\begin{rem} On a not\'e $M(-\lambda)^*$ le dual restreint de $M(-\lambda)$, le module de Verma de plus bas poids $\lambda$ (\cf par exemple le paragraphe qui pr\'ec\`ede la proposition 6 de \cite{Bry:OD}).
 \end{rem}

\begin{dem}
Choisissons une compactification r\'eguli\`ere $\G$ de $G$ et appliquons la proposition \ref{pro:sdc} . Aucun diviseur limitrophe ne contient l'orbite $BwB^-$. Ainsi, suivant les notations de la proposition pr\'ec\'edente :
$$I_{BwB^-} = \{p_{BwB^-}(D) \tq D \in {\cal D}\} \; , \;\; J_\Omega = \vide \p$$
De plus, d'apr\`es \cite[pro. A1, iv)]{B:Inf}, le r\'eseau engendr\'e par les caract\`eres $p_{B w B^-}(D)$,  $D \in {\cal D}$, est aussi le r\'eseau engendr\'e par les $T \croi T-$vecteurs propres de l'espace des fonctions rationnelles sur $\adh{T}$ (ou sur $T$). C'est donc : $\ma \{(\lambda ,-\lambda) \tq \lambda \in {\cal X}\}$. 
 \end{dem}

\tri

Dans les deux derni\`eres parties, on notera $U$ et $U^-$ les sous-groupes unipotents maximaux de $B$ et $B^-$. 

\section{Fin de la d\'emonstration du th\'eor\`eme principal }

D'apr\`es \cite[pro. A1, iv)]{B:Inf}, dans les compactifications de groupes r\'eductifs, toutes les $B \croi B^--$orbites sont de rang maximal.

On peut donc appliquer la proposition \ref{pro:sdc} dans ce cas.

Apr\`es avoir rappel\'e une param\'etrisation des $B \croi B^--$orbites de $\G$ (\ref{sec:cellu}), on va se sevir de la section pr\'ec\'edente pour \'eliminer une grande partie des termes du complexe de Grothendieck-Cousin (\cf la section \ref{ssec:tGC}) qui ne contribuent pas dans le calcul des multiplicit\'es des groupes de cohomologie $H^i(\G,\call{L}_h)$ (lemme \ref{lem:ww}). On exprimera alors ces multiplicit\'es \`a l'aide de groupes de cohomologie \`a support dans certaines $B \croi B^--$sous-vari\'et\'es de $\G$ (th\'eor\`eme \ref{thm:ctt}). Enfin, on pourra exprimer, \`a leur tour, ces groupes de cohomologie \`a support simplement en fonction de la restriction de $\call{L}_h$ \`a la vari\'et\'e torique $\adh{T} \cap \G_0$, d'\'eventail ${\cal E}^+$ (\ref{sec:ract}), et arriver ainsi \`a la formule du th\'eor\`eme \ref{thm:princi}.

\subsection{Les orbites dans la compactification}\label{sec:cellu}

Pour param\'etrer les $B \croi B^--$orbites de $\G$, on utilise les cellules et les $G \croi G-$orbites.

{\it On choisit une fois pour toutes un sous-groupe \`a un param\`etre $\zeta$ de $T \croi T$, dominant relativement \`a $B \croi B^-$ et tel que tous les points fix\'es par $\zeta$ le soient aussi par $T \croi T$.}

D'apr\`es \cite[pro. A1]{B:Inf}, les points fixes de $T \croi T$ sont dans les $G \croi G-$orbites ferm\'ees. Ce sont donc les points $(w,t).z_\sigma$ avec $(w,t) \in W \croi W$ et $\sigma$ un  c\^one  maximal  de  ${\cal E}^+$ (rappelons que $z_\sigma$ est le point-base associ\'e au c\^one $\sigma$ et que lorsque $\sigma$ d\'ecrit ${\cal E}^+$, l'ensemble des c\^ones maximaux, $z_\sigma$ d\'ecrit l'ensemble des points de $\G$ fix\'es par $B^- \croi B$).
Les $B \croi B^--$orbites de $\G$ sont d\`es lors les intersections non vides parmi les :

$${\cal O}_{w,t,\sigma} := {\cal O} \cap \G((w,t) . z_\sigma)$$
o\`u ${\cal O}$ est une $G \croi G-$orbite de $\G$, et $\G((w,t).z_\sigma)$ est la cellule de Bialynicki-Birula associ\'ee au point $(w,t).z_\sigma$ et au sous-groupe \`a un param\`etre $\zeta$ (\cf la section \ref{sec:vr}).

Pour la compactification magnifique $\adh{G_{ad}}$, il n'y a qu'un seul c\^one maximal : ${\cal C}^+$. On notera $\G_{ad}(w,t)$ la cellule de Bialynicki-Birula associ\'ee \`a $(w,t).z_{{\cal C}^+}$ et \`a $\zeta$.

On se servira plus loin de l'ouvert $\G_0$ form\'e des $x \in \G$ tels que l'orbite $B \croi B^- \hp x$ est ouverte dans $G \croi G \hp x$ (si $\G = \adh{G_{ad}}$, c'est la cellule ouverte). Plus g\'en\'eralement, on d\'efinit les relev\'es des cellules de $\adh{G_{ad}}$ dans $\G$ : suivant \cite[Pro. A2]{B:Inf}, comme $\G$ est r\'eguli\`ere, la surjection canonique $G \surj G_{ad}$ se prolonge en un morphisme $G \croi G-$\'equivariant $\pi : \G \to \adh{G_{ad}}$. Pour tous $w,t \in W$, soient $\G(w,t) : = \pi\inv (\G_{ad}(w,t))
$.

Posons aussi :
$$S(w,t) := \left( (w\inv,t\inv).\G(w,t) \right) \inter \adh{T} \cap \G_0$$
(si $w=t=1$, alors : $\G(1,1) = \G_0$ et $S(1,1) = \adh{T} \cap \G_0$).

Ces vari\'et\'es v\'erifient :
\begin{lem}[{\cite[\S 1.2 et 2.3]{BL:loc}}]\label{lem:loc}
Soient $w,t \in W$. Alors :\begin{itemize}
\item $\G(w,t)$ est une sous-vari\'et\'e ferm\'ee de $(w,t).\G_0$ ;
 \item $S(w,t)$ est une sous-vari\'et\'e ferm\'ee de $\adh{T} \cap \G_0$, form\'ee des $s \in \adh{T} \cap \G_0$ tels que la limite quand $a$ tend vers $0$ de :
$$(w\inv , t\inv)(\zeta)(a) . s$$
existe dans $\adh{T} \cap \G_0$ ;
\item on a la d\'ecomposition en $B \croi B^--$orbites :
$$\G(w,t) = \Dij_{ {\cal O} \: G\croi G-\mathrm{orbite} \atop \sigma \in {\cal E}^+ \: \mathrm{maximal}} {\cal O}_{w,t,\sigma} \;\; ;$$
\item l'application :
$$\appli{wU \croi tU^- \croi \adh{T} \cap \G_0}{ \G}{(a,b,s)}{(a,b).s}
$$
induit des isomorphismes de vari\'et\'es alg\'ebriques :
$$\begin{array}{ccccccc}
wU & \croi & tU^- & \croi & \adh{T} \cap \G_0 & \sta{\iso}{\to} & (w,t) \G_0 \\
\left( wU \cap Uw \right) & \croi & \left( t U^- \cap U^-t \right) & \croi & S(w,t) & \sta{\iso}{\to} & \G(w,t)
 \end{array} \; \Box$$
 \end{itemize}

 \end{lem}

\tri

Le dernier point du lemme a pour cons\'equence que, pour tout faisceau inversible et $\GXGT-$lin\'earis\'e sur $\G$, on a un isomorphisme de faisceaux $\ti{T} \croi \ti{T}-$lin\'earis\'es :
$$\call{L}_h\res{\G_0} \iso \fstr_U \cten \fstr_{U^-} \cten \call{L}_h\res{ \adh{T} \cap \G_0 } \p$$

\subsection{Calcul des multiplicit\'es des groupes de cohomologie \`a support dans les orbites}

Soit ${\cal O}_{w,t,\sigma}$ une $B \croi B^--$orbite de $\G$ de codimension $i$. On s'int\'eresse aux sous-quotients simples (comme $\goth g \croi \goth g-$modules) de dimension finie de $\ma H^i_{{\cal O}_{w,t,\sigma}}(\call{L}_h)$ et \`a leur multiplicit\'e. On d\'eduit de la proposition \ref{pro:sdc} p. \pageref{pro:sdc} le :

\begin{lem}\label{lem:ww}
Soit $L$ un $\GXGT-$module simple. Si la multiplicit\'e 
$$\left[ H^i_{{\cal O}_{w,t,\sigma}} (\call{L}_h) : L \right]$$
n'est pas nulle, alors : $w = t$ et il existe un caract\`ere $\mu \in \ti{\cal X}^+$  tel que $L = \End(L(\mu))$.
 \end{lem}

\begin{dem}
Soit $\chi$ le caract\`ere central du $\goth g \croi \goth g-$module $L$.

D'apr\`es la proposition \ref{pro:sdc}, si le $\goth g \croi \goth g-$module $\left( H^i_{{\cal O}_{w,t,\sigma}} \right)_\chi$ a une multiplicit\'e non nulle selon $L$, alors il en est de m\^eme pour un  $\goth g \croi \goth g-$module de la forme :
$$M^w(h_\sigma + \delta) \cten M^t(-h_\sigma - \delta)$$
o\`u $h_\sigma$ est le poids de la droite $\call{L}_h\res{z_\sigma}$ et, pour un certain diviseur $G \croi G-$invariant $D$ de $\G$, $\delta$ est celui de la droite ${\cal O}_X(D)\res{z_\sigma}$. Il est n\'ecessaire, aussi, que : 
$$(w,t) * (h_\sigma + \delta , -(h_\sigma + \delta)) \in W * \chi \p$$

Soit $(\lambda_1,\lambda_2)$ le plus haut poids de $L$ relativement \`a $B \croi B^-$ (\cad que $\lambda_1$ et $-\lambda_2$ sont dominants). On a alors :
$$w(h_\sigma + \delta + \rho) \in W(\lambda_1 + \rho) \et t(-h_\sigma-\delta-\rho) \in W  (-\lambda_2-\rho)$$
d'o\`u : $W(\lambda_1+ \rho) = W(\lambda_2+\rho)$ et : $\lambda_1 = -\lambda_2$ car $\lambda_1+\rho$ ainsi que $-(\lambda_2+ \rho)$ sont dominants et r\'eguliers. En cons\'equence, $L = \End(L(\lambda_1))$.

D'un autre c\^ot\'e, pour que le $\goth g \croi \goth g-$module $\ma M^w(h_\sigma + \delta) \cten M^t(-h_\sigma - \delta)$ ait un $\GXGT-$module parmi ses sous-quotients simples, il faut que $w(h_\sigma + \delta +\rho)$ et $t(h_\sigma + \delta +\rho)$ soient dominants r\'eguliers (\cf la remarque apr\`es le th\'eor\`eme \ref{thm:BoBry}). Cela entra\^{\i}ne que :
$$w(h_\sigma + \delta +\rho) = t(h_\sigma + \delta +\rho) =\lambda_1 = \lambda_2$$
et que : $w=t$.
\end{dem}

Gr\^ace au complexe de Grothendieck-Cousin, il r\'esulte imm\'ediatement de ce lemme que, pour tout $i \ge 0$, les repr\'esentations simples  de $\GXGT$ qui apparaissent dans la d\'ecomposition du $\GXGT-$module $H^i(\G , \call{L}_h )$ sont de la forme $\End(L(\mu))$ pour un certain $\mu \in \ti{\cal X}^+$. Ce lemme montre aussi que, dans le complexe de Grothendieck-Cousin et en vue du calcul des multiplicit\'es des modules $\End(L(\mu))$ dans les groupes de cohomologie $H^i(\G,\call{L}_h)$, on peut se passer des termes de la forme $H^i_{{\cal O}_{w,t,\sigma}} (\call{L}_h)$, avec $w \not=t$.  On se ram\`ene ainsi, dans la section suivante, \`a une \'etude des groupes de cohomologie \`a support dans les vari\'et\'es $\G(t,t)$, $t \in W$. Cette \'etude est plus facile que l'\'etude directe des groupes $H^i(\G,\call{L}_h)$ car, gr\^ace au lemme d'excision, on peut remplacer $\G$ par n'importe quel ouvert qui contient $\G(t,t)$.

\subsection{Simplification de l'expression des multiplicit\'es des groupes de cohomologie}

Avec les notations de la section pr\'ec\'edente (\ref{sec:cellu}), on a :
\begin{thm}\label{thm:ctt}
Pour tout $i \ge 0$ et tout $\mu \in \ti{\cal X}^+$ :
$$\left[ H^i(\G,\call{L}_h) : \End (L(\mu)) \right] = \sum_{t \in W} \left[ H^i_{\G(t,t)}(\call{L}_h) : \End(L(\mu)) \right] \p$$
 \end{thm}

\begin{dem}
Soit $L$ un $\GXGT-$module simple.

Suivant le th\'eor\`eme \ref{thm:cGC}, pour tout $i \ge 0$, $H^i(\G,\call{L}_h)$ est le $i-$\`eme groupe d'homologie du complexe de Grothendieck-Cousin :
$$K^* \; : \; 0 \to K^0 \sta{d^0}{\to} K^1 \sta{d^1}{\to} \sus \to K^{\dim \G} \to 0$$
o\`u pour chaque $p \ge 0$, $K^p$ est la somme directe :
$$\Plus_{{\cal O},w,t,\sigma} H^p_{{\cal O}_{w,t,\sigma}}(\call{L}_h)$$
index\'ee par les $B \croi B^--$orbites de $\G$ de codimension $p$.

\tri

Les applications $d^p$ d\'efinissent des morphismes (de $\goth g \croi \goth g-$modules) :
$$d^{\Omega,\Omega'} : H^p_\Omega(\call{L}_h) \to H^{p+1}_{\Omega'}(\call{L}_h)$$
pour toutes paires d'orbites ${\Omega,\Omega'}$ de codimensions $p$ et $p+1$.

Si, dans le complexe $K^*$, on ne garde que les termes de la forme $H_{{\cal O}_{t,t,\sigma}}(\call{L}_h)$ (d'apr\`es le lemme \ref{lem:ww}, ce sont les seuls qui puissent avoir une multiplicit\'e non nulle selon $L$), alors on obtient un nouveau complexe :
$$0 \to K^0_F \to K^1_F \to \sus$$
o\`u, pour tout $p \ge 0$ :
$$K^p_F := \Plus_{{\cal O},t,\sigma} H^p_{{\cal O}_{t,t,\sigma}}(\call{L}_h) \p$$
Les diff\'erentielles de ce complexe ont pour composantes les $d^{\Omega, \Omega'}$ o\`u $\Omega$ et ${\Omega'}$ sont des $B \croi B^--$orbites de la forme ${\cal O}_{t,t,\sigma}$.

De plus, on a :
$$\qq i \ge 0, \; [h^i(K^*):L] = [ h^i(K^*_F):L] \p$$

\tri

Pour des raisons de support, si un morphisme $d^{{\Omega, \Omega'}}$ est non nul, alors on a :
$$ {\Omega'} \sub \adh\Omega \et \codim({\Omega'} , \adh\Omega) = 1 \p$$

Or, gr\^ace \`a la description des adh\'erences des $B \croi B^--$orbites de $\G$ de \cite[th. du \S 2.1]{B:Inf} et \cite[pro. 2.4]{Sp:i}, on peut montrer le :

\begin{lem}[{\cite[pro. V.3.6]{these}}]
Soient ${\cal O,O'}$ des $G \croi G-$orbites de $\G$, $w,w' \in W$ et $\sigma , \sigma'$ des c\^ones maximaux de ${\cal E}^+$.

Alors, si :
$${\cal O'}_{w',w',\sigma'} \sub \adh{{\cal O}_{w,w,\sigma}} \et \codim({\cal O}_{w',w',\sigma'} , \adh{{\cal O}_{w,w,\sigma}}) =1 \;\;,$$
on a : $w = w'$ $\Box$
 \end{lem}

En cons\'equence, le complexe $K^*_F$ se d\'ecompose ainsi :
$$K^*_F = \Plus_{t \in W} K^*(t)$$
avec :
$$\qq p \ge 0, \, \qq t \in W, \;K^p(t) = \Plus_{{\cal O},\sigma} H^p_{{\cal O}_{t,t,\sigma}}(\call{L}_h) \;\; ;$$
il en r\'esulte aussi que :
$$\qq i \ge 0\,, \; [ h^i(K^*_F):L] = \sum_{t \in W}[h^i(K^*(t)) : L] \p$$

\tri 

Or, si $t \in W$, on reconna\^{\i}t en $K^*(t)$ le complexe de Grothendieck-Cousin associ\'e \`a $\call{L}_h$ et \`a la filtration :
$$\G \con \G(t,t) \con Z^0_t \con Z^1_t \con \sus $$
par les ferm\'es de $\G(t,t)$ d\'efinis par :
$$\qq i \ge 0 , \; Z^i_t := \Dij_{{\cal O},\sigma \atop \dim {\cal O}_{t,t,\sigma} \le \dim \G(t,t) -i} {\cal O}_{t,t,\sigma} \p$$

D\`es lors, d'apr\`es le th\'eor\`eme \ref{thm:cGC}, on a des isomorphismes de $\goth g \croi \goth g-$modules :
$$\qq i \ge 0 , \, \qq t \in W ,\: h^i(K^*(t)) \iso H^i_{\G(t,t)} (\call{L}_h) \p$$

\tri

En r\'esum\'e, si $i \ge 0$, alors :
$$[H^i(\G , \call{L}_h) : L ]  = [h^i(K^*) : L] = [ h^i(K^*_F) : L] = \sum_{t \in W} [h^i(K^*(t)) : L]$$
$$= \sum_{t \in W} [ H^i_{\G(t,t)} (\call{L}_h) : L] \p$$ 
 \end{dem}

\subsection{R\'eduction au cas torique}\label{sec:ract}

Pour conclure, il reste \`a d\'eterminer les multiplicit\'es 
$$\left[ H^i_{\G(t,t)}(\call{L}_h) : \End(L(\mu)) \right]$$
lorsque $t \in W$.

On va se ramener \`a au calcul du caract\`ere d'un groupe de cohomologie \`a support d'un faisceau inversible sur une vari\'et\'e torique.
 
Gr\^ace aux isomorphismes du lemme \ref{lem:loc}, on a un isomorphisme de $\ti{T} \croi \ti{T}-$modules :
$$H^i_{\G(t,t)} (\call{L}_h) \iso H^{l(t)}_{tU \cap Ut} (\fstr_{Ut}) \cten H^{l(t)}_{tU^- \cap U^-t} (\fstr_{U^-t}) \tens H^{i - 2l(t)}_{S(t,t)} (\call{L}_h \res{\adh{T} \cap \G_0}) \p $$

Et, en utilisant les isomorphimes de $\ti{T}-$modules suivants (\cf \cite[lem. 12.8]{Ke:GCc} ou \cite[lemme 3.16]{Ku:BGG} pour les deux premiers) :
$$H^{l(t)}_{tUt\inv \cap U} (\fstr_{tUt\inv}) \iso M(t\rho - \rho)$$
(le module de Verma de plus haut poids $t\rho -\rho$)
$$H^{l(t)}_{tU^-t\inv \cap U^-} (\fstr_{tU^-t\inv}) \iso M(-t\rho + \rho)$$
(le module de Verma de plus bas poids $-t\rho + \rho$), et
$$H^{i-2l(t)}_{S(t,t)} (\call{L}_h\res{\adh{T}\cap \G_0} ) \iso \Plus_{\nu \in \ti{\cal X}} m_\nu \C_\nu$$
(pour certains entiers $m_\nu$ \infra{
On consid\`ere $H^{i-2l(t)}_{S(t,t)} (\call{L}_h\res{\adh{T}\cap \G_0} )$ comme un $\ti{T}-$module gr\^ace \`a l'action de $\ti{T} \croi \{1\}$ sur $\adh{T} \cap \G_0$ et sur $S(t,t)$.}
), on obtient un isomorphisme de $\ti{T} \croi \ti{T}-$modules :
$$H^i_{\G(t,t)} (\call{L}_h) \iso \Plus_{\nu \in \ti{\cal X}} m_\nu \Bigm( M(t(\nu + \rho) -\rho) \cten M(-(t(\nu + \rho) -\rho)) \Bigm) \p$$

Mais, comme les caract\`eres des $\goth g \croi \goth g-$modules simples $L$ (\`a plus hauts poids) sont $\Z-$lin\'eairement ind\'ependants, on trouve que, pour tout caract\`ere $\mu \in \ti{\cal X}^+$, la multiplicit\'es selon $\End(L(\mu))$ de $H^i_{\G(t,t)} (\call{L}_h)$ v\'erifie :
$$\left[ H^i_{\G(t,t)} (\call{L}_h) \tq \End(L(\mu)) \right] = \sum_{ \nu \in W * \mu} m_\nu \left[ M(t(\nu+ \rho) - \rho) \tq L(\mu) \right]$$
$$ = m_{t\inv (\mu + \rho ) - \rho} $$
(\cf la remarque qui suit le th\'eor\`eme \ref{thm:BoBry}).

D\`es lors, pour tout $i \ge 0$ et pour tout $\mu \in \ti{\cal X}^+$, on a :
$$\left[ H^i(\G , \call{L}_h) \tq \End(L(\mu) \right] = m_\mu + \sum_{t\in W , \, t \not=1 } m_{t\inv * \mu} \p$$

\tri

Enfin, on v\'erifie que les vari\'et\'es toriques $\adh{T} \cap \G_0$ et $\adh{T} \cap \G_0 \moins S(t,t)$ ont pour \'eventails :
${\cal E}^+ \mbox{ avec } | {\cal E}^+ | = {\cal C}^+$
et un ensemble ${\cal E}^+_t$ avec $\ma |{\cal E}^+_t| = \uni_{\alpha \in J_{t^{-1}} } \alpha^\ort \cap {\cal C}^+$. 

La formule de \cite[th. 2.6]{O:cbag} (\cf aussi \cite[th. II.4.2]{these}) permet alors d'exprimer les entiers $m_{t\inv * \mu}$ en fonction de $h$ et des sous-espaces topologiques ${\cal C}^+$ et $ \uni_{\alpha \in J_t} \alpha^\ort \cap {\cal C}^+ $ de ${\cal Y_\R}$.

Cela ach\`eve la d\'emonstration du th\'eor\`eme \ref{thm:princi}. 
\tri
En fait, pour la cohomologie des fibr\'es en droites sur les vari\'et\'es toriques, la m\'ethode du complexe de Grothendieck-Cousin aboutit directement \`a une formule qui fait intervenir la cohomologie d'Ishida (\cf \cite[th. II.4.4]{these}).

\vskip 1cm

{\it Je remercie mon directeur de th\`ese, Michel Brion, pour sa relecture attentive et patiente, et pour les am\'eliorations et corrections qu'il m'a fait apporter aux premi\`eres versions. Un grand merci \'egalement \`a Syu Kato pour les bons \'echanges que nous avons eus.}



\end{document}